\documentclass[a4paper,12pt,reqno]{amsart}
\usepackage{amsmath}
\usepackage{amssymb}
\usepackage{amsthm}

\usepackage{amscd}
\usepackage{amsfonts}
\usepackage{setspace}
\usepackage{version}

\title[On finding many solutions to $S$-unit equations]{On finding many solutions to $S$-unit equations by solving linear equations on average}
\author{Adam J Harper}
\address{Department of Pure Mathematics and Mathematical Statistics, Wilberforce Road, Cambridge CB\textup{3} \textup{0}WA, England}
\email{A.J.Harper@dpmms.cam.ac.uk}
\date{18th August 2011}
\thanks{The author is supported by a studentship from the Engineering and Physical Sciences Research Council of the United Kingdom.}

\numberwithin{equation}{section}
\theoremstyle{plain}

\onehalfspace
\newcommand{\Z}{\mathbb{Z}}

\newtheorem{prop1}{Proposition}
\newtheorem{thm1}{Theorem}
\newtheorem{thm2}[thm1]{Theorem}
\newtheorem{lem1}{Lemma}
\newtheorem{lem2}[lem1]{Lemma}
\newtheorem{nthres}{Number Theory Result}
\newtheorem{nthres2}[nthres]{Number Theory Result}
\newtheorem{nthres3}[nthres]{Number Theory Result}
\newtheorem{nthres4}[nthres]{Number Theory Result}
\newtheorem{thm3}[thm1]{Theorem}
\newtheorem{thm4}[thm1]{Theorem}

\begin{document}

\maketitle

\begin{abstract}
We give improved lower bounds for the number of solutions of some $S$-unit equations over the integers, by counting the solutions of some associated linear equations as the coefficients in those equations vary over sparse sets. This method is quite conceptually straightforward, although its successful implementation involves, amongst other things, a slightly subtle use of a large sieve inequality. We also present two other results about solving linear equations on average over their coefficients.
\end{abstract}

\section{Introduction}
In this paper we will give bounds for the number of solutions of some $S$-unit equations over the integers. That is, we will investigate the integer solutions of equations such as
$$ a + 1 = c, $$
where all prime factors of the variables are required to lie in some given finite set $\mathcal{S}$. (One can equally well speak of $S$-unit equations in any algebraic number field, with $\mathcal{S}$ being a set of primes of the field, but we shall not consider the more general situation.) Throughout we will use the letter $s$ to denote the cardinality of the relevant set $\mathcal{S}$.

$S$-unit equations have been studied fairly intensively for at least thirty years. For example, in 1984 Evertse~\cite{ev} used Diophantine approximation methods to give upper bounds, depending only on $s$ (and, in the number field case, also the degree of the field), for the number of $S$-unit solutions to equations like $a+1=c$. In 1988, Erd\H{o}s, Stewart and Tijdeman~\cite{erdsttijd} used methods from combinatorial number theory, in particular some results they proved about the largest prime factor of a product of sums, to show the existence of sets $\mathcal{S}$ for which $a+b=c$ has many $S$-unit solutions with $a,b,c$ {\em coprime} integers. In 2003, Evertse, Moree, Stewart and Tijdeman~\cite{emst} gave an extension of Erd\H{o}s, Stewart and Tijdeman's result to weighted $S$-unit equations in $n$ variables\footnote{When comparing Evertse, Moree, Stewart and Tijdeman's paper~\cite{emst} with other results, one should carefully note that they allow the variables in their equations to take {\em rational}, and not just integer, values. Their results about the equation $a_{1}+a_{2}+...+a_{n}=1$, $n \geq 2$, with rational variables $a_{i}$, translate into results about the integer equation $a_{1}+a_{2}+...+a_{n}=a_{n+1}$.} (as well as some other results). Many more upper bound results are now known, and there has also been much work in other directions.

In their 2007 paper~\cite{konsound}, Konyagin and Soundararajan used a nice and simple combinatorial argument to improve the lower bound result of Erd\H{o}s, Stewart and Tijdeman~\cite{erdsttijd}, and a more involved argument, requiring lower bounds for the quantity of very smooth numbers\footnote{A number is said to be $y$-{\em smooth} if all of its prime factors are at most $y$. By ``very smooth'' we mean ``$y$-smooth, with $y$ very small relative to other quantities in the argument''.} in certain arithmetic progressions, to show the existence of sets $\mathcal{S}$ for which $a+1=c$ has many $S$-unit solutions. Finally we mention two recent papers~\cite{lagsound1,lagsound2} of Lagarias and Soundararajan, in which the Generalised Riemann Hypothesis is used (powering a version of the circle method) to give lower bounds for the number of coprime $S$-unit solutions to $a+b=c$, in the special case where $\mathcal{S}$ is fixed as the set of the first $s$ prime numbers.

In this paper we will give lower bounds for the number of solutions of some $S$-unit equations, using general results about character sums and Kloosterman sums.

\vspace{12pt}
Before stating our results, we explain a general strategy for exhibiting many solutions to an $S$-unit equation. The work of Erd\H{o}s, Stewart and Tijdeman~\cite{erdsttijd}, of Evertse, Moree, Stewart and Tijdeman~\cite{emst}, and of Konyagin and Soundararajan~\cite{konsound} can be seen to fit into this strategy, and our approach is also based upon it.

It seems a natural guess that, for given $s$, an $S$-unit equation such as $a+b=c$ will have most coprime solutions (at least roughly) if $\mathcal{S}$ consists of the first $s$ primes. For example, one might imagine that the fine distribution of $S$-units looks fairly random, subject to global constraints on their density, and so one is most likely to find solutions when the $S$-units are as ``concentrated'' as possible (i.e. when the elements of $\mathcal{S}$ are small). It is reasoning of this kind that gives rise to many conjectures about the number of solutions to $S$-unit equations. More precisely, it is known that for any fixed $a \geq 1$,
$$ \#\{n \leq x : \textrm{if prime } p \mid n \textrm{ then } p \leq \log^{a}x \} = x^{1 - 1/a + o(1)} \;\;\;\;\; \textrm{ as } x \rightarrow \infty. $$ 
See $\S 2$ for a proof of a lower bound of this kind. If we choose $a = 3/2 + \epsilon$, and let $x$ be large, then\footnote{As is usual, we employ the notation $\Omega(\epsilon)$ to denote any function of $\epsilon$ that is $\gg \epsilon$.} the right hand side is $x^{1/3 + \Omega(\epsilon) + o(1)}$, and we might imagine that
\begin{eqnarray}
&& \#\{(a,b,c) \in \Z^{3} : \gcd(a,b,c)=1; a+b=c; \textrm{ if prime } p \mid abc \textrm{ then } p \leq \log^{3/2+\epsilon}x \} \nonumber \\
& \approx & x^{1/3 + \Omega(\epsilon) + o(1)} \cdot x^{1/3 + \Omega(\epsilon) + o(1)} \cdot x^{1/3 + \Omega(\epsilon) + o(1)} / x = x^{\Omega(\epsilon) + o(1)}, \nonumber
\end{eqnarray}
which is $\geq e^{(\log^{3/2+\epsilon}x)^{2/3 - \epsilon}}$. The conjectured lower bound for the maximum number of coprime $S$-unit solutions to $a+b=c$ is, accordingly, $e^{s^{2/3-\epsilon}}$. See e.g. Konyagin and Soundararajan's paper~\cite{konsound} for further calculations of this kind.

A standard approach for counting solutions to equations whose variables lie in complicated sets is the circle method. However, it is well known that this almost never works for equations with only two free variables, and for equations with three variables one would usually work with sets having density at least $\log^{-C}x$, say, for $C$ some positive constant. In contrast, we would want to work with sets having density more like $x^{-1/2}$ or $x^{-2/3}$, and about which we also don't have much Fourier-analytic information. In their papers~\cite{lagsound1,lagsound2}, Lagarias and Soundararajan assume the Generalised Riemann Hypothesis to obtain suitable information.

We will avoid these difficulties by giving up on the idea of counting solutions for $\mathcal{S}$ that is fixed in advance, and instead allow ourselves to vary $\mathcal{S}$ a little in response to things that happen in our arguments. In other words, we will start out with a set $\mathcal{S}'$, and will find lots of ``almost solutions'' to the equation we are interested in whose prime factors come from $\mathcal{S}'$. Then we will add a few extra primes into $\mathcal{S}'$, forming a set $\mathcal{S}$, in such a way that quite a few of the almost solutions are turned into actual solutions with prime factors from $\mathcal{S}$. Thus when Erd\H{o}s, Stewart and Tijdeman~\cite{erdsttijd} were seeking many coprime solutions to $a + b = c$, they took $\mathcal{S}'$ to consist of all primes up to some point, and used their results about prime factors of sums to show the existence of a fixed (and not too large) $b$ such that $a + b = c$ has many solutions in $a,c$ having all their prime factors in $\mathcal{S}'$. Then they set $\mathcal{S} = \mathcal{S}' \cup \{p : p \textrm{ prime}, p \mid b\}$. Konyagin and Soundararajan~\cite{konsound} studied the $S$-unit equation $a + b = c$ by studying the solutions of
$$ c - b = au, $$
where (roughly speaking) $a,b,c$ had only small prime factors (from $\mathcal{S}'$), and $u$ was any small integer. They showed the existence of so many solutions that some ``popular'' $u$ had to occur in many of them, and then set $\mathcal{S} = \mathcal{S}' \cup \{p : p \textrm{ prime}, p \mid u\}$. Similarly, they obtained $S$-unit solutions of $a + 1 =c$ by studying solutions of
$$ c - 1 = au. $$

We will study the $S$-unit equation $a + 1 = c$ by studying the solutions, in smooth numbers $a,c$ and integers $u,w$ from certain intervals (roughly speaking), of the equation
$$ au + 1 = cw. $$
Similarly we will study the $S$-unit equation $a + b + 1= c$ by studying the solutions of $au + b + 1 = cw$, with $a,b,c$ smooth numbers and $u,w$ integers of small modulus.

\vspace{12pt}
As a warm-up, in $\S 2$ we shall prove the following result:
\begin{prop1}[Erd\H{o}s, Stewart and Tijdeman, 1988]
There are arbitrarily large sets $\mathcal{S}$ of prime numbers such that
$$ \#\{(a,b,c) \in \Z^{3} : \gcd(a,b,c)=1; a+b=c; \textrm{ if prime } p \mid abc \textrm{ then } p \in \mathcal{S}\} \geq e^{(1/2\sqrt{2})\sqrt{s/\log s}},  $$
where $s := \#\mathcal{S}$.
\end{prop1}
Actually Erd\H{o}s, Stewart and Tijdeman~\cite{erdsttijd} proved this with the constant $1/2\sqrt{2}$ replaced by $4-o(1)$, as $s \rightarrow \infty$, and Konyagin and Soundararajan~\cite{konsound} showed a stronger lower bound $e^{s^{2-\sqrt{2}-o(1)}}$. See below for more about this. However, our proof is a little different from existing proofs, and motivates the proofs of our subsequent results. In $\S 2$ we shall also briefly discuss the $n$-variable $S$-unit equation $a_{1}+a_{2}+...+a_{n}=0$.

Our main theorem concerns the two variable $S$-unit equation that we have already mentioned several times.
\begin{thm1}
Let $\epsilon > 0$ be fixed. There are arbitrarily large sets $\mathcal{S}$ of prime numbers such that
$$ \#\{(a,c) \in \Z^{2} : a + 1 = c, \textrm{ and if prime } p \mid ac \textrm{ then } p \in \mathcal{S}\} \geq e^{s^{1/6 - \epsilon}}, $$
where $s := \#\mathcal{S}$.
\end{thm1}
This will be proved in $\S 3$, and improves the previous best result of Konyagin and Soundararajan~\cite{konsound}, who obtained this with the exponent $1/6 - \epsilon$ replaced by $1/16$ (or in fact a number around $1/15.88$). We will also show that if a widely believed conjecture about short character sums is true then one could prove this result with an exponent $1/5 - \epsilon$, and with an easier proof. The conjectured lower bound here is $e^{s^{1/2 - \epsilon}}$.

In $\S 4$, we shall prove:
\begin{thm2}
Let $\epsilon > 0$ be fixed. There are arbitrarily large sets $\mathcal{S}$ of prime numbers such that
$$ \#\{(a,b,c) \in \Z^{3} : a + b + 1 = c; \textrm{ if prime } p \mid abc \textrm{ then } p \in \mathcal{S}; a,b,-c \neq -1\} \geq e^{s^{\lambda_{0} - \epsilon}}, $$
where $\lambda_{0} \approx 0.53551$ is the real root of the cubic $4x^{3}-5x^{2}+9x-4$, and $s := \#\mathcal{S}$.
\end{thm2}
We specify that our solutions will have $a,b,-c \neq -1$ to rule out large classes of ``degenerate'' solutions, such as $a - a + 1 = 1$. The conjectured lower bound here is $e^{s^{2/3-\epsilon}}$, whereas e.g. Erd\H{o}s, Stewart and Tijdeman's argument~\cite{erdsttijd} for the equation $a+b=c$ can be directly adapted to give a lower bound $e^{s^{1/2-\epsilon}}$. Konyagin and Soundararajan's work~\cite{konsound} on the equation $a+b=c$, which the author was unable to improve in that case, doesn't seem to extend to handle the equation $a+b+1=c$, because shifting by $1$ (or any other non-zero quantity) seems to spoil one of their applications of the pigeonhole principle. As before, assuming a certain conjecture about short character sums allows an improved exponent, that is $\approx 0.55496$, in Theorem 2.

\vspace{12pt}
As indicated in our earlier discussion, our proofs involve estimating the number of small integer solutions to certain linear equations as the coefficients in an equation vary over ``arithmetically interesting'' sets. This is a topic that has already received quite a lot of attention: see e.g. the recent paper of Shparlinski~\cite{shpar}, and references therein, as well as the section on linear equations in Shparlinski's survey paper~\cite{shpar2}.

The results in $\S\S 3-4$ are somewhat tailored to the proofs of Theorems 1 and 2, in that various parameters are optimised for those applications, but the techniques of the proofs should transfer to other contexts and may be of independent interest. In $\S 5$ we present two general results on solving linear equations on average, which extend some results in the literature but don't seem to lead to sharper conclusions about $S$-unit equations. We also sketch an argument that, in a very limited sense, the arguments of $\S 3$ are somewhat optimal for obtaining conclusions about the $S$-unit equation $a+1=c$, which rules out a few obvious approaches to extending Theorem 1.

\section{Proof of Proposition 1}
To prove Proposition 1 we shall require two fairly easy lemmas:
\begin{lem1}[Part of Siegel's lemma, 1929]
Let $n \geq 2$, and let $\alpha_{1},...,\alpha_{n}$ be integers bounded in absolute value by $B \geq 1$. Then there exist integers $z_{1},...,z_{n}$, not all zero and bounded in absolute value by $(nB)^{1/(n-1)}$, such that
$$ \alpha_{1}z_{1} + \alpha_{2}z_{2} + ... + \alpha_{n}z_{n} = 0. $$
\end{lem1}

Lemma 1 is a consequence of the pigeonhole principle: we follow the proof given in Chapter 2.4 of Baker's book~\cite{baker} to obtain the clean bound $(nB)^{1/(n-1)}$. For any $C \geq 1$, as $y_{1},...,y_{n}$ vary over non-negative integers bounded by $C$ the linear form $\sum_{i=1}^{n} \alpha_{i}y_{i}$ takes at most $n[CB]+1$ distinct values, noting that each term $\alpha_{i}y_{i}$ grows the set of possible values in one direction only, according to the sign of $\alpha_{i}$. Here $[\cdot]$ denotes integer part. Thus if
$$ ([C]+1)^{n} > n[CB]+1 $$
we must have $\sum_{i=1}^{n} \alpha_{i}y_{i} = \sum_{i=1}^{n} \alpha_{i}y_{i}'$ for some $(y_{1},...,y_{n}) \neq (y_{1}',...,y_{n}')$, and then we can set $(z_{1},...,z_{n}) = (y_{1}-y_{1}',...,y_{n}-y_{n}')$. Since $([C]+1)^{n} > [C]([C]+1)^{n-1} + 1$ and $n[CB]+ 1 \leq nCB+1$, we see the desired inequality is satisfied when $C = [(nB)^{1/(n-1)}]$.

The full version of Siegel's lemma supplies small solutions to systems of linear equations: see Baker's book~\cite{baker}, which also has much other interesting material.

\begin{lem2}
Let $\nu > 0$. If $x$ is sufficiently large, depending on $\nu$, then the following is true.
If $1 \leq a \leq 100$ (say), and $\mathcal{T} \subseteq [2,\log^{a}x]$ is any set of prime numbers, and if $b \geq \nu$ is such that $\#\mathcal{T} = (\log^{b+1}x)/(a\log\log x)$, then
$$ \#\{ n \leq x : n \textrm{ is squarefree, with all of its prime factors from } \mathcal{T} \} \geq \frac{x^{b/a + 1/(2a\log\log x)}}{6 \log^{b+1} x}. $$
\end{lem2}

If we set $k = [(\log x)/(a \log\log x)]$, it will clearly suffice to prove such a lower bound for $\binom{\#\mathcal{T}}{k}$. But using Stirling's formula, since $x$ is assumed to be large we have
$$ \binom{\#\mathcal{T}}{k} \geq (\#\mathcal{T}/k)^{k} \frac{\#\mathcal{T}^{\#\mathcal{T}-k}}{3(\#\mathcal{T}-k)^{\#\mathcal{T}-k} \sqrt{k}} \geq (\#\mathcal{T}/k)^{k} \frac{e^{k/2}}{3 \sqrt{k}}. $$
Then one calculates that
$$ (\#\mathcal{T}/k)^{k} \geq ((\#\mathcal{T} a \log\log x)/\log x)^{(\log x)/(a \log\log x)} (\#\mathcal{T}/k)^{-1} = x^{b/a} (k/\#\mathcal{T}), $$
and the inequality in Lemma 2 follows with some room to spare (using the trivial bounds $k/\sqrt{k} = \sqrt{k} \geq 1$ and $\#\mathcal{T} \leq \log^{b+1}x$).

\vspace{12pt}
Now we can swiftly deduce Proposition 1. By e.g. the prime number theorem, if $x$ is a large parameter we can find $3$ disjoint sets $\mathcal{T}_{1},\mathcal{T}_{2},\mathcal{T}_{3}$ of prime numbers satisfying
$$ \mathcal{T}_{i} \subseteq [2,\log^{2 + (2\log 4)/\log\log x}x] \;\; \textrm{ and } \;\; \#\mathcal{T}_{i} \geq \frac{\log^{2 + (2\log 4)/\log\log x}x}{4 \log(\log^{2 + (2\log 4)/\log\log x}x)} = \frac{\log^{b_{i}+1}x}{\log(\log^{2 + (2\log 4)/\log\log x}x)}, $$
where $b_{i} = 1 + (\log 4)/\log\log x$. In view of Lemma 2, the subsets $\mathcal{A}_{i}$ of $[1,x]$ consisting of squarefree integers with all their prime factors from $\mathcal{T}_{i}$ each have size at least $x^{1/2 + 1/(4\log\log x) + O(1/(\log\log x)^{2})}$.

Using Lemma 1, we see that
\begin{eqnarray}
\sum_{\alpha_{1} \in \mathcal{A}_{1}} \sum_{\alpha_{2} \in \mathcal{A}_{2}} \sum_{\alpha_{3} \in \mathcal{A}_{3}} \sum_{|z_{1}| \leq \sqrt{3x},...,|z_{3}| \leq \sqrt{3x}, \atop (z_{1},z_{2},z_{3}) \neq (0,0,0)} \textbf{1}_{\alpha_{1}z_{1}+\alpha_{2}z_{2}+\alpha_{3}z_{3} = 0} & \geq & \sum_{\alpha_{1} \in \mathcal{A}_{1}} \sum_{\alpha_{2} \in \mathcal{A}_{2}} \sum_{\alpha_{3} \in \mathcal{A}_{3}} 1 \nonumber \\
& \geq & x^{3/2 + 3/(4\log\log x) + O(1/(\log\log x)^{2})}. \nonumber
\end{eqnarray}
In fact this is still true if one only sums over tuples $(z_{1},z_{2},z_{3})$ with no coordinate equal to zero: tuples with two or more zero coordinates make no contribution anyway, and e.g. the equation $\alpha_{2}z_{2} + \alpha_{3}z_{3} = 0$ can have at most one solution with $\alpha_{2} \in \mathcal{A}_{2}, \alpha_{3} \in \mathcal{A}_{3}$ and $z_{2},z_{3} \neq 0$ fixed, since such $\alpha_{2}$ and $\alpha_{3}$ are always coprime (by construction). Thus the total contribution from such ``degenerate'' tuples is at most $O(x(\#\mathcal{A}_{1} + \#\mathcal{A}_{2} + \#\mathcal{A}_{3})) = o((\#\mathcal{A}_{1})(\#\mathcal{A}_{2})(\#\mathcal{A}_{3}))$.

If we swap the order of the summations on the left hand side, this implies that for some popular $(z_{1},z_{2},z_{3})$ (with no zero coordinates) we must have
$$ \sum_{\alpha_{1} \in \mathcal{A}_{1}} \sum_{\alpha_{2} \in \mathcal{A}_{2}} \sum_{\alpha_{3} \in \mathcal{A}_{3}} \textbf{1}_{\alpha_{1}z_{1}+...+\alpha_{3}z_{3} = 0} \geq \frac{x^{3/2 + 3/(4\log\log x) + O(1/(\log\log x)^{2})}}{(2\sqrt{3x})^{3}} \geq x^{3/(4\log\log x) + O(1/(\log\log x)^{2})}. $$
Now we shall fix this popular tuple $(z_{1},z_{2},z_{3})$. The numbers $\alpha_{1}z_{1},...,\alpha_{3}z_{3}$ may sometimes have a divisor in common, but we know this isn't the case for any $\alpha_{1},...,\alpha_{3}$, so if we reduce each tuple $(\alpha_{1}z_{1},...,\alpha_{3}z_{3})$ by dividing out the highest common factor we will obtain distinct triples of coprime integers.

At this point we have found at least $e^{3\log x/(4\log\log x) + O(\log x/(\log\log x)^{2})}$ triples of coprime integers solving the equation $A+B+C=0$, and having all their prime factors from the set
$$ \{p \textrm{ prime} : p \leq \log^{2 + (2\log 4)/\log\log x}x\} \cup \{p \textrm{ prime} : p|z_{1}z_{2}z_{3} \}. $$
This has size at most
$$ \frac{1.1 \log^{2 + (2\log 4)/\log\log x}x}{\log(\log^{2 + (2\log 4)/\log\log x}x)} + 3\log(\sqrt{3x})/\log 2 \leq \frac{17.6 \log^{2}x}{2\log\log x} + 3\log(\sqrt{3x})/\log 2 \leq \frac{8.9\log^{2}x}{\log\log x} $$
when $x$ is large, and (setting $s$ equal to this, and comparing with $e^{3\log x/(4\log\log x)}$) our claimed lower bound $e^{(1/2\sqrt{2})\sqrt{s/\log s}}$ immediately follows.
\begin{flushright}
Q.E.D.
\end{flushright}

For general $n \geq 3$, Evertse, Moree, Stewart and Tijdeman~\cite{emst} showed (if we translate from their statements about rational variables to statements about integer equations) that there exist arbitrarily large sets $\mathcal{S}$ for which
$$ a_{1}+a_{2}+...+a_{n} = 0 $$
has at least $e^{s^{1-1/(n-1)-o(1)}}$ non-degenerate $S$-unit solutions. In fact they had a bit more precise bound. Here non-degenerate means not only that the variables have highest common factor one, but also that no proper subsum of the variables vanishes. This is implied by the highest common factor condition when $n=3$, but not for larger $n$, where one wants to rule out large numbers of uninteresting solutions such as $a-a+1-1=0$. By the reasoning explained in the introduction, it is conjectured that one should be able to have $e^{s^{1-1/n-o(1)}}$ non-degnerate solutions, so (in a certain sense) our knowledge is not too far from the presumed truth when $n$ is large.

Our proof could be adapted to give a lower bound $e^{s^{1-1/(n-1)-o(1)}}$ for all $n$, by considering $n$ disjoint sets $\mathcal{T}_{i} \subseteq [2,\log^{(n-1)/(n-2) + \epsilon}x]$ of primes, but it is not obvious how it could be made to respect the condition of no vanishing subsums. In view of Evertse, Moree, Stewart and Tijdeman's work~\cite{emst}, it probably isn't much worth pursuing this.

\section{Proof of Theorem 1}

\subsection{Preliminary observations}
Let $\delta > 0$ be a parameter, which we shall take to be very small. Also let $W,X,Z$ be numbers (that should be thought of as large) satisfying $W \leq \min\{X,Z\}$ and $X^{1/100} \leq Z \leq X^{100}$, say. Following the general strategy set out in the introduction, we shall prove Theorem 1 by showing that, if $\mathcal{C} \subseteq [X^{1-\delta},X]$ and $\mathcal{A} \subseteq [Z^{1-\delta},Z]$ are sets having certain properties,
$$ \sum_{a \in \mathcal{A}} \sum_{c \in \mathcal{C}} \sum_{1 \leq w \leq W} \textbf{1}_{cw-au=1 \textrm{ for some } u} \geq \frac{X^{1+101\delta} W^{2}}{Z}. $$
Note that if $cw-au=1$ then $1 \leq u = (cw-1)/a \leq XW/Z^{1-\delta}$, so if the above holds then for some popular $1 \leq w \leq W$ and $u$ we must have
$$ \sum_{a \in \mathcal{A}} \sum_{c \in \mathcal{C}} \textbf{1}_{cw-au=1} \geq \frac{X^{101\delta}}{Z^{\delta}} \geq X^{\delta}. $$
We will show that, in particular, one can choose $\mathcal{A}, \mathcal{C}$ to consist of integers having all their prime factors from a set $\mathcal{S}'$ of size about $\log^{6+\epsilon}X$, where $\epsilon > 0$ is the quantity from the statement of Theorem 1. We will then set
$$ \mathcal{S} := \mathcal{S}' \cup \{p : p \textrm{ prime}, p \mid uw\}, $$
so that $\#\mathcal{S} \leq \#\mathcal{S}' + (\log W)/\log 2 + (\log(XW/Z^{1-\delta}))/\log 2 \leq \log^{6+2\epsilon} X$ if $X$ is large enough, and
\begin{eqnarray}
&& \#\{(A,C) \in \Z^{2} : A + 1 = C, \textrm{ and if prime } p \mid AC \textrm{ then } p \in \mathcal{S}\} \nonumber \\
& \geq & \#\{(a,c) \in \Z^{2} : au + 1 = cw, \textrm{ and if prime } p \mid ac \textrm{ then } p \in \mathcal{S}'\} \nonumber \\
& \geq & e^{\delta(\#\mathcal{S})^{1/(6+2\epsilon)}} \nonumber \\
& \geq & e^{(\#\mathcal{S})^{1/6 - \epsilon}}, \nonumber
\end{eqnarray}
as claimed in Theorem 1. Note that we can make the set $\mathcal{S}$ arbitrarily large by choosing $W,X,Z$ to be correspondingly large.

To obtain the desired lower bound on the triple sum we shall use a discrete form of the circle method, i.e. we shall decompose the indicator function $\textbf{1}_{cw-au=1 \textrm{ for some } u}$ as a sum of ``harmonics\footnote{In the language of the very nice book of Iwaniec and Kowalski~\cite{iwkow}.}''. It turns out that decomposing into Dirichlet characters is more effective than decomposing into additive characters: thus we have
\begin{eqnarray}
\sum_{a \in \mathcal{A}} \sum_{c \in \mathcal{C}} \sum_{1 \leq w \leq W} \textbf{1}_{cw-au=1 \textrm{ for some } u} & = & \sum_{a \in \mathcal{A}} \sum_{c \in \mathcal{C}} \sum_{1 \leq w \leq W}  \frac{1}{\phi(a)} \sum_{\chi \textrm{ mod } a} \chi(c) \chi(w) \overline{\chi}(1) \nonumber \\
& = & \sum_{a \in \mathcal{A}} \frac{1}{\phi(a)} \#\{c \in \mathcal{C} : (c,a) = 1\} \#\{1 \leq w \leq W : (w,a)=1\} + \nonumber \\
&& + O(\sum_{a \in \mathcal{A}} \frac{1}{\phi(a)} \sum_{\chi \textrm{ mod } a, \atop \chi \neq \chi_{0}} \left|\sum_{c \in \mathcal{C}} \chi(c) \right| \left|\sum_{1 \leq w \leq W} \chi(w) \right|), \nonumber
\end{eqnarray}
where $\chi_{0}$ denotes the principal Dirichlet character to modulus $a$. We will show that, under suitable conditions, the ``big Oh'' term is negligible compared with the first term, and also the first term is at least (twice) as big as the claimed lower bound $X^{1+101\delta}W^{2}/Z$.

Before we launch into estimating the ``big Oh'' term, we observe that (as is hopefully suggested by our array of notation) it is likely to be beneficial to take $\mathcal{A}$ and $\mathcal{C}$ on rather different ``scales'', i.e. to choose $Z$ and $X$ to have different sizes. The author found this quite counterintuitive at first, since when solving the equation $A+1=C$ one clearly must take $A$ and $C$ to be about the same size. However, given the presence of the extra variables $u,w$ it turns out that some asymmetry can be very helpful: indeed, if we could show (as we might at best hope to) that $\left|\sum_{c \in \mathcal{C}} \chi(c) \right| \ll_{\eta} (\#\mathcal{C})^{1/2+\eta}$ for most non-principal characters $\chi$, for any $\eta > 0$ (with the implicit constant depending on $\eta$), we would clearly have saved most if $\mathcal{C}$ was very large compared with $a$. Also in some of our arguments it will be very helpful if we can ``factor'' the set $\mathcal{C}$ as a product, to give ourselves flexibility in arranging the sums we must bound. This will be explained more precisely in the following sections.

Versions of the circle method involving Dirichlet characters are no doubt employed in many papers, but the author would particularly like to reference the work of Garaev and coauthors: the reader might consult e.g. a paper of Cilleruelo and Garaev~\cite{cillgar} on the least number with totient in an arithmetic progression, and a recent preprint of Garaev~\cite{garaev}, amongst much other work. These papers also emphasise the importance of working with variables on appropriate scales.

\subsection{A conditional proof}
First we shall give a fairly straightforward proof of Theorem 1, but conditional on the far from proven conjecture that if $\eta > 0$, $\chi$ is any non-principal character mod $a$, and $W$ is arbitrary, then
$$ \left| \sum_{1 \leq w \leq W} \chi(w) \right| \ll_{\eta} \sqrt{W} a^{\eta}. $$
Actually, assuming this conjecture to hold for all $a$ we shall prove the theorem with the stronger lower bound $e^{s^{1/5-\epsilon}}$.

Assume that we have sets $\mathcal{Q} \subseteq [Q^{1-\delta},Q]$ and $\mathcal{R} \subseteq [R^{1-\delta},R]$, where $QR=X$, such that all of the products $qr$, $q \in \mathcal{Q}$, $r \in \mathcal{R}$ are distinct, and are exactly the elements of the set $\mathcal{C} \subseteq [X^{1-\delta},X]$. Then for any $a \in \mathcal{A}$,
\begin{eqnarray}
\frac{1}{\phi(a)} \sum_{\chi \textrm{ mod } a, \atop \chi \neq \chi_{0}} \left|\sum_{c \in \mathcal{C}} \chi(c) \right| \left|\sum_{1 \leq w \leq W} \chi(w) \right| & = & \frac{1}{\phi(a)} \sum_{\chi \textrm{ mod } a, \atop \chi \neq \chi_{0}} \left|\sum_{q \in \mathcal{Q}} \chi(q) \right| \left|\sum_{r \in \mathcal{R}} \chi(r) \right| \left|\sum_{1 \leq w \leq W} \chi(w) \right| \nonumber \\
& \ll_{\eta} & \frac{\sqrt{W} a^{\eta}}{\phi(a)} \sum_{\chi \textrm{ mod } a, \atop \chi \neq \chi_{0}} \left|\sum_{q \in \mathcal{Q}} \chi(q) \right| \left|\sum_{r \in \mathcal{R}} \chi(r) \right| \nonumber \\
& \leq & \sqrt{W} a^{\eta} \sqrt{ \frac{1}{\phi(a)} \sum_{\chi \textrm{ mod } a} \left|\sum_{q \in \mathcal{Q}} \chi(q) \right|^{2} \cdot \frac{1}{\phi(a)} \sum_{\chi \textrm{ mod } a} \left|\sum_{r \in \mathcal{R}} \chi(r) \right|^{2} }, \nonumber
\end{eqnarray}
using the conjecture and the Cauchy--Schwarz inequality. Now if $Q$ and $R$ are smaller than $Z^{1-\delta}$, and therefore (recall $\S 3.1$) smaller than $a$, then the elements of $\mathcal{Q}$ will all be distinct modulo $a$, similarly for the elements of $\mathcal{R}$. Therefore if we expand the squares in the last line, and just perform the summations over $\chi$, we will obtain a bound
$$ \ll_{\eta} \sqrt{W} a^{\eta} \sqrt{ \#\mathcal{Q} \cdot \#\mathcal{R} } = \sqrt{W} a^{\eta} \sqrt{\#\mathcal{C}} \leq Z^{\eta} \sqrt{W} \sqrt{\#\mathcal{C}}, $$
which is the best we could have hoped to establish.

\vspace{12pt}
Before we go on, we shall introduce some useful non-standard notation. We will write $B(\delta)$ to denote any quantity such that, for any $\sigma > 0$, we have $X^{-\sigma} \leq B(\delta) \leq X^{\sigma}$ when $X$ is large enough and $\delta$ is small enough in terms of $\sigma$. From place to place $B(\delta)$ will denote different such quantities.

With this notation, it remains to determine how we can set the relative sizes of $W,X$ and $Z$, and construct the sets $\mathcal{A}$ and $\mathcal{C}$, so that both
$$ \sum_{a \in \mathcal{A}} \frac{1}{\phi(a)} \#\{c \in \mathcal{C} : (c,a) = 1\} \#\{1 \leq w \leq W : (w,a)=1\} \geq B(\delta) \#\mathcal{A} \sqrt{W} \sqrt{\#\mathcal{C}}, $$
(which is at least $2 \sum_{a \in \mathcal{A}} (1/\phi(a)) \sum_{\chi \textrm{ mod } a, \chi \neq \chi_{0}} |\sum_{c \in \mathcal{C}} \chi(c)| |\sum_{1 \leq w \leq W} \chi(w)|$, by the foregoing calculations), and also
$$ \sum_{a \in \mathcal{A}} \frac{1}{\phi(a)} \#\{c \in \mathcal{C} : (c,a) = 1\} \#\{1 \leq w \leq W : (w,a)=1\} \geq \frac{B(\delta) X W^{2}}{Z}. $$
We shall choose $\mathcal{C}$ in such a way that all of its elements are coprime to all elements of $\mathcal{A}$. Moreover we shall certainly end up having $W \geq Z^{1/1000}$, say, so by the sieve of Eratosthenes--Legendre (as in e.g. chapter 3.1 of Montgomery and Vaughan~\cite{mv}) we will have
$$ \#\{1 \leq w \leq W : (w,a)=1\} = W \prod_{p \mid a} (1-1/p) + O(d(a)) \geq (1/2)W \prod_{p \mid a} (1-1/p) $$
(if $Z$ is large enough, and $d(\cdot)$ denotes the divisor function). Thus it will suffice to have
$$ \sqrt{W} \sqrt{\#\mathcal{C}} \geq B(\delta) Z \;\;\; \textrm{ and } \;\;\; \#\mathcal{A} \#\mathcal{C} \geq B(\delta) XW. $$
Finally if $\alpha \geq 0$ is any number such that $\#\mathcal{A} \geq Z^{1-\alpha}$ and $\#\mathcal{C} \geq X^{1-\alpha}$ both hold, then {\em it will suffice to have}
$$ W^{1/2}X^{(1-\alpha)/2} \geq B(\delta) Z \;\;\; \textrm{ and } \;\;\; Z^{1-\alpha} X^{-\alpha} \geq B(\delta) W. $$

Comparing the two inequalities that we wish to satisfy, we must have
$$ Z \leq B(\delta) X^{(1-\alpha)/2} (Z^{1-\alpha}X^{-\alpha})^{1/2} \Rightarrow Z^{(1+\alpha)/2} \leq B(\delta) X^{1/2-\alpha} \Rightarrow Z \leq B(\delta) X^{(1-2\alpha)/(1+\alpha)}. $$
However, earlier in our argument we supposed that $X = QR \leq (Z^{1-\delta})^{2} = B(\delta) Z^{2}$, so the inequality can only possibly be satisfied (for large $X$ and small $\delta$) if
$$ (1-2\alpha)/(1+\alpha) \geq 1/2, $$
i.e. if $\alpha \leq 1/5$. On the other hand, if we have $\alpha = 1/5 - \epsilon$; and $Z = X^{-\epsilon}X^{(1-2\alpha)/(1+\alpha)}$; and $W = X^{-\epsilon}X^{(1-4\alpha+\alpha^{2})/(1+\alpha)}$; and $Q=R=X^{1/2}$; and if $\delta$ is sufficiently small in terms of $\epsilon$, and $X$ is sufficiently large; then we will have
$$ W^{1/2}X^{(1-\alpha)/2} = X^{-\epsilon/2}X^{(1-2\alpha)/(1+\alpha)} \geq X^{-\epsilon} B(\delta) X^{(1-2\alpha)/(1+\alpha)} = B(\delta) Z, $$
and similarly $Z^{1-\alpha} X^{-\alpha} \geq B(\delta) W$ and $Q,R \leq Z^{1-\delta}$, as we wanted.

Finally we can let $\mathcal{T}_{1}, \mathcal{T}_{2}, \mathcal{T}_{3}$ be three disjoint sets of primes satisfying
$$ \mathcal{T}_{i} \subseteq [(\log^{5+\epsilon}X)/2,\log^{5+\epsilon}X] \;\;\; \textrm{ and } \;\;\; \#\mathcal{T}_{i} \geq \frac{\log^{5+\epsilon}X}{10 \log(\log^{5+\epsilon}X)}, $$
and let $\mathcal{Q} \subseteq [1,Q]$, $\mathcal{R} \subseteq [1,R]$ and $\mathcal{A} \subseteq [1,Z]$ be the subsets of squarefree numbers having all their prime factors from $\mathcal{T}_{1},\mathcal{T}_{2},\mathcal{T}_{3}$ respectively. By Lemma 2, these sets have size at least $Q^{1-1/5+\Omega(\epsilon)+o(1)}, R^{1-1/5+\Omega(\epsilon)+o(1)}, Z^{1-1/5+\Omega(\epsilon)+o(1)}$, and one can check in the proof of Lemma 2 that this will still hold if we remove numbers smaller than $Q^{1-\delta},R^{1-\delta},Z^{1-\delta}$ (respectively) from the sets. Now we can just define
$$ \mathcal{C} = \{qr : q \in \mathcal{Q}, r \in \mathcal{R}\}, $$
where clearly all the products $qr$ are distinct. This completes our construction of sets $\mathcal{A}, \mathcal{C}$ producing many solutions of $cw-au=1$, and all of whose elements have their prime factors from a set of size $\leq \log^{5+\epsilon}X$.

\subsection{An unconditional proof}
Now we shall prove Theorem 1 without the aid of any conjecture. To do this we shall require the following two results:
\begin{nthres}[P\'{o}lya--Vinogradov inequality, P\'{o}lya, Vinogradov, 1918]
Let $\chi$ be a non-principal character to modulus $q$, having conductor $r$. For any integers $M$ and $N$ with $N > 0$,
$$ \left| \sum_{n=M+1}^{M+N} \chi(n) \right| \leq d(q/r) \sqrt{r} \log r, $$
where $d(\cdot)$ denotes the divisor function.
\end{nthres}

\begin{nthres2}[Multiplicative Large Sieve inequality, Bombieri, 1965]
Let $Y < Z$ be integers, and let $\mathcal{Q}$ be any finite set of positive integers. Then for any complex numbers $a_{n}$,
$$ \sum_{q \in \mathcal{Q}} \frac{1}{\phi(q)} \sum_{\chi \textrm{ mod } q} |\tau(\chi)|^{2} \left|\sum_{Y < n \leq Z} \chi(n) a_{n} \right|^{2} \leq 7 D(\mathcal{Q}) \max\{Z-Y, M(\mathcal{Q})^{2}\} \sum_{Y < n \leq Z} d(n) |a_{n}|^{2}, $$
where $d(\cdot)$ denotes the divisor function, $D(\mathcal{Q}) := \max_{q \in \mathcal{Q}} d(q)$, $M(\mathcal{Q}) := \max_{q \in \mathcal{Q}} q$, and $\tau(\chi)$ denotes the Gauss sum of $\chi$.
\end{nthres2}

Probably both of these results are fairly familiar to the reader. Most books seem to state the P\'{o}lya--Vinogradov inequality without reference to the conductor of $\chi$, claiming a weaker bound $\ll \sqrt{q} \log q$ for the character sum. However, most (perhaps all?) proofs of the P\'{o}lya--Vinogradov inequality really prove the stronger statement: see e.g. chapter 9.4 of Montgomery and Vaughan~\cite{mv}. Number Theory Result 2 is proved in Bombieri's paper~\cite{bomb}, but here too most books seem to state a restricted version where the summations are over primitive characters $\chi$ only. Perhaps surprisingly, it will prove quite important to have access to Number Theory Results 1 and 2 in the less usual forms that we have stated.  

We proceed again to bound the ``big Oh'' term from the end of $\S 3.1$. We still assume that we can factor the set $\mathcal{C}$ as a product of $\mathcal{Q}$ and $\mathcal{R}$, as in $\S 3.2$, so that
\begin{eqnarray}
\sum_{a \in \mathcal{A}} \frac{1}{\phi(a)} \sum_{\chi \textrm{ mod } a, \atop \chi \neq \chi_{0}} \left|\sum_{c \in \mathcal{C}} \chi(c) \right| \left|\sum_{1 \leq w \leq W} \chi(w) \right| & = & \sum_{a \in \mathcal{A}} \frac{1}{\phi(a)} \sum_{\chi \textrm{ mod } a, \atop \chi \neq \chi_{0}} \left|\sum_{q \in \mathcal{Q}} \chi(q) \right| \left|\sum_{r \in \mathcal{R}} \chi(r) \right| \left|\sum_{1 \leq w \leq W} \chi(w) \right| \nonumber \\
& \leq & \sqrt{\sum_{a \in \mathcal{A}} \frac{1}{\phi(a)} \sum_{\chi \textrm{ mod } a} \left|\sum_{q \in \mathcal{Q}} \chi(q) \right|^{2}} \cdot \nonumber \\
&& \cdot \sqrt{\sum_{a \in \mathcal{A}} \frac{1}{\phi(a)} \sum_{\chi \textrm{ mod } a, \atop \chi \neq \chi_{0}} \left|\sum_{r \in \mathcal{R}} \chi(r) \right|^{2} \left|\sum_{1 \leq w \leq W} \chi(w) \right|^{2}}. \nonumber
\end{eqnarray}
If we assume, as before, that $\mathcal{Q} \subseteq [Q^{1-\delta},Q]$ for some $Q \leq Z^{1-\delta}$, and we also apply Number Theory Result 1 together with a crude bound $\ll_{\eta} Z^{\eta}$ for the divisor function, we find the right hand side is
$$ \ll_{\eta} Z^{\eta} \sqrt{\#\mathcal{A} \#\mathcal{Q}} \sqrt{\sum_{a \in \mathcal{A}} \frac{1}{\phi(a)} \sum_{\chi \textrm{ mod } a, \atop \chi \neq \chi_{0}} \textrm{cond}(\chi) \left|\sum_{r \in \mathcal{R}} \chi(r) \right|^{2} }. $$
Here we wrote $\textrm{cond}(\chi)$ to denote the conductor of the character $\chi$. However, if the modulus $a$ is squarefree, which we will arrange to be the case for all $a \in \mathcal{A}$, then we just have $\textrm{cond}(\chi) = |\tau(\chi)|^{2}$ (see e.g. chapter 9.2 of Montgomery and Vaughan~\cite{mv}). This means that we can use Number Theory Result 2 to estimate the above, obtaining a bound
$$ \ll_{\eta} Z^{2\eta} \sqrt{\#\mathcal{A} \#\mathcal{Q}} \sqrt{\max\{R,Z^{2}\} \#\mathcal{R}} = Z^{2\eta} \sqrt{\#\mathcal{A} \#\mathcal{C}} \max\{\sqrt{R},Z\}. $$

\vspace{12pt}
Exactly similarly to $\S 3.2$, we now need to arrange that
$$ \sum_{a \in \mathcal{A}} \frac{1}{\phi(a)} \#\{c \in \mathcal{C} : (c,a) = 1\} \#\{1 \leq w \leq W : (w,a)=1\} \geq B(\delta) \sqrt{\#\mathcal{A} \#\mathcal{C}} \max\{\sqrt{R},Z\} $$
and that
$$ \sum_{a \in \mathcal{A}} \frac{1}{\phi(a)} \#\{c \in \mathcal{C} : (c,a) = 1\} \#\{1 \leq w \leq W : (w,a)=1\} \geq \frac{B(\delta)X W^{2}}{Z}. $$
Working as in $\S 3.2$, and letting $\alpha \geq 0$ be such that  $\#\mathcal{A} \geq Z^{1-\alpha}$ and $\#\mathcal{C} \geq X^{1-\alpha}$, as there, it will suffice if
$$ X^{(1-\alpha)/2} W \geq B(\delta) Z^{(1+\alpha)/2} \max\{\sqrt{R},Z\} \;\;\; \textrm{ and } \;\;\; Z^{1-\alpha} X^{-\alpha} \geq B(\delta) W. $$
Comparing the inequalities, if they are to be satisfied we must have
$$ Z^{(1+\alpha)/2} \max\{\sqrt{R},Z\} \leq B(\delta) X^{(1-\alpha)/2}(Z^{1-\alpha}X^{-\alpha}) \Rightarrow Z^{(3\alpha-1)/2} \max\{\sqrt{R},Z\} \leq B(\delta) X^{(1-3\alpha)/2}; $$
and if we impose the constraint that $R \leq Z^{2}$ this forces that $Z \leq B(\delta) X^{(1-3\alpha)/(1+3\alpha)}$. On the other hand, we supposed earlier that $Q \leq Z^{1-\delta}$, and so $X=QR \leq B(\delta) Z^{3}$. This means that, under our constraint $R \leq Z^{2}$, the inequalities can only possibly be satisfied if
$$ (1-3\alpha)/(1+3\alpha) \geq 1/3, $$
i.e. if $\alpha \leq 1/6$. The reader may check that if we tried to set $R > Z^{2}$, then we would end up requiring that
$$ Z^{(3\alpha-1)/2} \sqrt{R} \leq B(\delta) (QR)^{(1-3\alpha)/2} \leq B(\delta) Z^{(1-3\alpha)/2} R^{(1-3\alpha)/2} \Rightarrow R \leq B(\delta) Z^{2(1-3\alpha)/3\alpha}, $$  
which still cannot be satisfied for any value of $\alpha$ larger than $1/6$.

Finally we observe that setting $\alpha = 1/6 - \epsilon$; and $Z = X^{-\epsilon} X^{(1-3\alpha)/(1+3\alpha)}$; and $W = X^{-\epsilon} X^{(1-5\alpha)/(1+3\alpha)}$; and $Q=Z^{1-\delta}$ (and $R=X/Q < Z^{2}$); our various inequalities and assumptions {\em will} all be satisfied. Theorem 1 follows on constructing the sets $\mathcal{Q}, \mathcal{R}, \mathcal{A}$ exactly similarly to $\S 3.2$, simply replacing $5+\epsilon$ by $6+\epsilon$ in all of the exponents.
\begin{flushright}
Q.E.D.
\end{flushright}

\section{Proof of Theorem 2}

\subsection{More preliminary observations}
Our proof of Theorem 2 will follow much the same strategy as our proof of Theorem 1, so we shall be brief in recalling this. This time we shall have a small parameter $\delta > 0$, and large quantities $W,X,Y,Z$ satisfying $Z \leq \min\{X,Y\}$ and $X^{1/100} \leq Y \leq X^{100}$, say. We shall prove Theorem 2 by showing that, if $\mathcal{C} \subseteq [X^{1-\delta},X]$, $\mathcal{B} \subseteq [Y^{1-\delta},Y]$ and $\mathcal{A} \subseteq [Z^{1-\delta},Z]$ are sets having certain properties,
$$ \sum_{a \in \mathcal{A}} \sum_{b \in \mathcal{B}} \sum_{c \in \mathcal{C}} \sum_{1 \leq w \leq W} \textbf{1}_{cw-au=b+1 \textrm{ for some } u \neq 0} \geq X^{\delta}W \left(\frac{Y+XW}{Z^{1-\delta}} \right). $$

Note that we only want to count quadruples $(a,b,c,w)$ such that $cw-au=b+1$ for some {\em non-zero} $u$, so that distinct values of $a$ give rise to distinct values of $au$. In $\S 3$ we didn't need to specify this, since when $cw-au=1$ it was necessarily the case that $u \neq 0$. We will be able to force that $u \neq 0$ without any effort, since it will transpire that choosing $Y=(XW)^{1/(1-\delta)}=B(\delta)XW$ is roughly optimal, (recall the notation $B(\delta)$ from $\S 3.2$), and then
$$ |cw| \leq XW = Y^{1-\delta} < b + 1, $$
so if $cw-au=b+1$ then $u$ must be non-zero. We will not fix $Y$ in this way yet, so that the reader can see the optimality of the choice in our calculations, but we will assume that $Y \leq B(\delta) XW$. This is a simplifying assumption that should not hurt us, since for given $Y,Z,W$ and $X < B(\delta)Y/W$ we would expect to only increase our quadruple sum, and not essentially change the quantity on the right hand side, by increasing $X$. 

In any event, if $cw-au=b+1$ then $-Y/Z^{1-\delta} \leq u = (cw-b-1)/a \leq XW/Z^{1-\delta}$, so if the above holds then for some popular $1 \leq w \leq W$ and $u \neq 0$ we must have
$$ \sum_{a \in \mathcal{A}} \sum_{b \in \mathcal{B}} \sum_{c \in \mathcal{C}} \textbf{1}_{cw-au=b+1} \geq X^{\delta}. $$
We will show that, in particular, one can choose $\mathcal{A}, \mathcal{B}, \mathcal{C}$ to consist of integers having all their prime factors from a set $\mathcal{S}'$ of size $\log^{(1/\lambda_{0})+\epsilon}X$, where $\lambda_{0}$ and $\epsilon > 0$ are the quantities from the statement of Theorem 2. We will then set $\mathcal{S} := \mathcal{S}' \cup \{p : p \textrm{ prime}, p \mid uw\}$, and Theorem 2 will follow exactly as Theorem 1 did in $\S 3.1$.

Once again we shall use Dirichlet characters to implement the circle method, noting
\begin{eqnarray}
&& \sum_{a \in \mathcal{A}} \sum_{b \in \mathcal{B}} \sum_{c \in \mathcal{C}} \sum_{1 \leq w \leq W} \textbf{1}_{cw-au=b+1 \textrm{ for some } u} \nonumber \\
& \geq & \sum_{a \in \mathcal{A}} \sum_{b \in \mathcal{B}, \atop (a,b+1)=1} \sum_{c \in \mathcal{C}} \sum_{1 \leq w \leq W} \textbf{1}_{cw-au=b+1 \textrm{ for some } u} \nonumber \\
& = & \sum_{a \in \mathcal{A}} \frac{1}{\phi(a)} \#\{b \in \mathcal{B} : (b+1,a)=1\} \#\{c \in \mathcal{C} : (c,a) = 1\} \#\{1 \leq w \leq W : (w,a)=1\} + \nonumber \\
&& + O(\sum_{a \in \mathcal{A}} \frac{1}{\phi(a)} \sum_{\chi \textrm{ mod } a, \atop \chi \neq \chi_{0}} \left|\sum_{b \in \mathcal{B}} \overline{\chi}(b+1) \right| \left|\sum_{c \in \mathcal{C}} \chi(c) \right| \left|\sum_{1 \leq w \leq W} \chi(w) \right|), \nonumber
\end{eqnarray}
where $\chi_{0}$ denotes the principal Dirichlet character to modulus $a$. We will show that, under suitable conditions, the ``big Oh'' term is negligible compared with the first term, and also the first term is at least (twice) as big as our claimed lower bound.

As in $\S 3$ we will be working with variables on different scales, and indeed $X,Y,Z$ will all end up being different to one another. However, since we are already summing over a third variable, namely $b \in \mathcal{B}$, we will not ``factor'' any of our sets as we did in the proof of Theorem 1. There is one extra observation that we need to make, which is crucial to the proof of Theorem 2. We may assume that for each integer $n \neq 0$,
$$ \#\{(c,c') \in \mathcal{C}^{2} : c-c'=n\} \ll_{\eta} X^{\eta} \;\;\; \textrm{ and } \;\;\; \#\{(c,c',c'',c''') \in \mathcal{C}^{4} : cc'-c''c'''=n\} \ll_{\eta} X^{\eta} $$
for any $\eta > 0$, similarly for elements from the set $\mathcal{B}$. For if the first bound failed for some $n \neq 0$ and some $\eta > 0$, then since $c-c'+1=n+1$ we would already have $\gg X^{\eta}$ non-degnerate solutions $(c,-c',n+1)$ to our target equation. The same is true if the second bound failed, just noting that at most $O_{\eta}(X^{\eta/2})$ of the quadruples $(c,c',c'',c''')$ can produce the same products $cc', c''c'''$, using a crude bound for the divisor function. (And we will be choosing $\mathcal{C},\mathcal{B}$ to consist of integers with all their prime factors from small sets, so these solutions would indeed be the many $S$-unit solutions that we wanted).

\subsection{A conditional proof}
Similarly to $\S 3.2$, we shall first prove Theorem 2 on the conjecture that if $\eta > 0$, $\chi$ is any non-principal character mod $a$, and $W$ is arbitrary, then
$$ \left| \sum_{1 \leq w \leq W} \chi(w) \right| \ll_{\eta} \sqrt{W} a^{\eta}. $$
Actually we shall prove the theorem with the stronger lower bound $e^{s^{\lambda_{1}-\epsilon}}$ by assuming this conjecture, where $\lambda_{1} \approx 0.55496$ is the root of the cubic $x^{3}-2x^{2}-x+1$ that lies between one half and one.

For any $a \in \mathcal{A}$ we have
\begin{eqnarray}
\frac{1}{\phi(a)} \sum_{\chi \textrm{ mod } a, \atop \chi \neq \chi_{0}} \left|\sum_{b \in \mathcal{B}} \overline{\chi}(b+1) \right| \left|\sum_{c \in \mathcal{C}} \chi(c) \right| \left|\sum_{1 \leq w \leq W} \chi(w) \right| & \ll_{\eta} & \frac{\sqrt{W} a^{\eta}}{\phi(a)} \sum_{\chi \textrm{ mod } a, \atop \chi \neq \chi_{0}} \left|\sum_{b \in \mathcal{B}} \chi(b+1) \right| \left|\sum_{c \in \mathcal{C}} \chi(c) \right| \nonumber \\
& \leq & \sqrt{W} a^{\eta} \sqrt{\frac{1}{\phi(a)} \sum_{\chi \textrm{ mod } a} \left|\sum_{b \in \mathcal{B}} \chi(b+1) \right|^{2}} \cdot \nonumber \\
&& \cdot \sqrt{\frac{1}{\phi(a)} \sum_{\chi \textrm{ mod } a} \left|\sum_{c \in \mathcal{C}} \chi(c) \right|^{2}}, \nonumber
\end{eqnarray}
using the conjecture and the Cauchy--Schwarz inequality. Expanding the squares and performing the sums over characters, we find this is at most
$$ \sqrt{W} a^{\eta} \sqrt{\left( \#\mathcal{B} + \sum_{b \in \mathcal{B}} \sum_{b' \in \mathcal{B}, b' \neq b} \textbf{1}_{a | b-b'} \right) \left(\#\mathcal{C} + \sum_{c \in \mathcal{C}} \sum_{c' \in \mathcal{C}, c' \neq c} \textbf{1}_{a | c-c'} \right)}. $$
Since all of the numbers $b-b'$ are non-zero integers between $-Y$ and $Y$, and we may assume (recall $\S 4.1$) that no such integer has more than $O_{\eta}(Y^{\eta})$ representations as $b-b'$, the double sum over $b,b'$ is $\ll_{\eta} Y^{\eta} (1+Y/a) \ll Y^{1+\eta}/Z^{1-\delta}$, and similarly the sum over $c,c'$ is $\ll_{\eta} X^{1+\eta}/Z^{1-\delta}$. Here we used our assumption, from $\S 4.1$, that $Z \leq \min\{X,Y\}$.

\vspace{12pt}
Recalling our useful notation $B(\delta)$ from $\S 3.2$, and our assumption from $\S 4.1$ that $Y \leq B(\delta) XW$, it remains to determine how we can set the relative sizes of $W,X,Y$ and $Z$, and construct the sets $\mathcal{A},\mathcal{B},\mathcal{C}$, so that both
\begin{eqnarray}
&& \sum_{a \in \mathcal{A}} \frac{1}{\phi(a)} \#\{b \in \mathcal{B} : (b+1,a)=1\} \#\{c \in \mathcal{C} : (c,a) = 1\} \#\{1 \leq w \leq W : (w,a)=1\} \nonumber \\
& \geq & B(\delta) \sqrt{W} \#\mathcal{A} \sqrt{(\#\mathcal{B} + Y/Z)(\#\mathcal{C} + X/Z)}, \nonumber
\end{eqnarray}
(which is at least $2 \sum_{a \in \mathcal{A}} (1/\phi(a)) \sum_{\chi \textrm{ mod } a, \chi \neq \chi_{0}} |\sum_{b \in \mathcal{B}} \overline{\chi}(b+1)| |\sum_{c \in \mathcal{C}} \chi(c)| |\sum_{1 \leq w \leq W} \chi(w)|$, by the foregoing calculations), and also
\begin{eqnarray}
&& \sum_{a \in \mathcal{A}} \frac{1}{\phi(a)} \#\{b \in \mathcal{B} : (b+1,a)=1\} \#\{c \in \mathcal{C} : (c,a) = 1\} \#\{1 \leq w \leq W : (w,a)=1\} \nonumber \\
& \geq & \frac{B(\delta) XW^{2}}{Z}. \nonumber
\end{eqnarray}
We shall choose $\mathcal{C}$ in such a way that all of its elements are coprime to all elements of $\mathcal{A}$, and we shall choose $\mathcal{B}$ such that $\#\{b \in \mathcal{B} : (b+1,a)=1\} \geq \#\mathcal{B}/2$ for each $a \in \mathcal{A}$. Thus, exactly similarly to $\S 3.2$, it will suffice to have
$$ \sqrt{W} \#\mathcal{B} \#\mathcal{C} \geq B(\delta) Z \sqrt{(\#\mathcal{B} + Y/Z)(\#\mathcal{C} + X/Z)} \;\;\; \textrm{ and } \;\;\; \#\mathcal{A} \#\mathcal{B} \#\mathcal{C} \geq B(\delta) XW. $$
Finally if $\alpha \geq 0$ is any number such that $\#\mathcal{A} \geq Z^{1-\alpha}$, $\#\mathcal{B} \geq Y^{1-\alpha}$ and $\#\mathcal{C} \geq X^{1-\alpha}$ all hold, then {\em it will suffice to have}
$$ \sqrt{W} Y^{1-\alpha} X^{1-\alpha} \geq B(\delta) Z \sqrt{(Y^{1-\alpha} + Y/Z)(X^{1-\alpha} + X/Z)} \;\;\; \textrm{ and } \;\;\; Z^{1-\alpha}Y^{1-\alpha} \geq B(\delta) X^{\alpha}W. $$

It will be easier to think about the first constraint if we replace it by the four constraints that come from multiplying out the brackets, namely
$$ \sqrt{W} Y^{1-\alpha} X^{1-\alpha} \geq B(\delta) Z \sqrt{ \max\{Y^{1-\alpha}X^{1-\alpha},XY^{1-\alpha}/Z,X^{1-\alpha}Y/Z,XY/Z^{2}\} }. $$
It is clearly best to choose $W$ as large as possible subject to the other constraint, i.e. to take $W = B(\delta)X^{-\alpha}Y^{1-\alpha}Z^{1-\alpha}$; this leaves us to satisfy
$$ Y^{3(1-\alpha)/2} X^{1-(3\alpha/2)} \geq B(\delta) Z^{(1+\alpha)/2} \sqrt{ \max\{Y^{1-\alpha}X^{1-\alpha},XY^{1-\alpha}/Z,X^{1-\alpha}Y/Z,XY/Z^{2}\} }. $$
Now if we examine the powers of $Y$ on the left and right hand sides, on the left side we have a power $3(1-\alpha)/2$, and on the right the largest power of $Y$ is $1/2$. Since we shall certainly end up having $\alpha \leq 2/3$, (and in fact we shall end up with $\alpha = \lambda_{1}-\epsilon$, of course), the power on the left hand side will be larger, so it will be best if we choose $Y$ as large as possible. We have been supposing since $\S 4.1$ that $Y \leq B(\delta) XW$, so we should take $Y=B(\delta) XW=B(\delta)X^{1-\alpha}Y^{1-\alpha}Z^{1-\alpha}$, i.e. $Y^{\alpha}=B(\delta)X^{1-\alpha} Z^{1-\alpha}$. Noting that $X^{1-\alpha}Y/Z \geq XY^{1-\alpha}/Z$ when we make this choice (since $Y \geq X$), so that one of the four terms in our maximum becomes redundant, we are left to satisfy
\begin{eqnarray}
Z^{3(1-\alpha)^{2}} X^{2\alpha-3\alpha^{2}+3(1-\alpha)^{2}} & \geq & B(\delta) Z^{\alpha(1+\alpha)}  \max\{Z^{(1-\alpha)^{2}}X^{1-\alpha},X^{(\alpha+1)(1-\alpha)}Z^{1-2\alpha},XZ^{1-3\alpha}\} \nonumber \\
& = & B(\delta) \max\{Z^{1-\alpha+2\alpha^{2}}X^{1-\alpha},X^{(\alpha+1)(1-\alpha)}Z^{1-\alpha+\alpha^{2}},XZ^{(1-\alpha)^{2}}\}. \nonumber
\end{eqnarray}
(Here we raised both sides of our previous inequality to the power $2\alpha$, to clear denominators in the exponents).

On shifting all the powers of $X$ and $Z$ to the right hand side, this remaining constraint becomes
$$ \max\{Z^{-2+5\alpha-\alpha^{2}}X^{-2+3\alpha},X^{-2+4\alpha-\alpha^{2}}Z^{-2+5\alpha-2\alpha^{2}},X^{-2+4\alpha}Z^{-2(1-\alpha)^{2}}\} \leq B(\delta). $$
The powers of $Z$ in the first two terms are non-negative for all $1/2 \leq \alpha \leq 1$, so it will be best to choose $Z$ as small as possible such that the third term is smaller than $B(\delta)$, i.e. to choose $Z = B(\delta) X^{(2\alpha-1)/(1-\alpha)^{2}}$. Making this choice leaves us to satisfy
$$ \max\{X^{(2\alpha-1)(-2+5\alpha-\alpha^{2}) + (-2+3\alpha)(1-\alpha)^{2}},X^{(-2+4\alpha-\alpha^{2})(1-\alpha)^{2} + (2\alpha-1)(-2+5\alpha-2\alpha^{2})}\} \leq B(\delta), $$
where we raised both sides of our previous inequality to the power $(1-\alpha)^{2}$. Since $X$ is supposed to be a large parameter, this just amounts to asking for the maximum of the two exponents to be negative, i.e. for
$$ \max\{-2\alpha+3\alpha^{2}+\alpha^{3},-\alpha+\alpha^{2}+2\alpha^{3}-\alpha^{4}\} < 0, $$
where $1/2 \leq \alpha \leq 1$. Cancelling one multiple of $\alpha$, we want that $\alpha^{2}+3\alpha-2 < 0$ and that $\alpha^{3}-2\alpha^{2}-\alpha+1 > 0$. The reader may check that the second of these inequalities is a stronger condition (i.e. forces $\alpha$ to be smaller), and it will be satisfied when $\alpha = \lambda_{1} - \epsilon$, by definition of $\lambda_{1}$ as the relevant root of that cubic. If $\alpha$ is chosen in this way; and $Z,Y,W$ are chosen in terms of $X$ in the way listed above, that is if
$$ Z = B(\delta) X^{(2\alpha-1)/(1-\alpha)^{2}} \;\;\; \textrm{ and } \;\;\; Y = B(\delta) X^{\alpha/(1-\alpha)} \;\;\; \textrm{ and } \;\;\; W = B(\delta) X^{(2\alpha-1)/(1-\alpha)}, $$
for suitable quantities $B(\delta)$; and if $X$ is sufficiently large in terms of $\epsilon$, and $\delta$ is sufficiently small; then the reader may check back that all of our constraints (including that $Z \leq \min\{X,Y\}$) will be satisfied.

Finally, and similarly to $\S 3.2$, we can let $\mathcal{T}_{2}$ be the set of all primes strictly smaller than $(\log^{(1/\lambda_{1})+\epsilon}X)/2$, and let $\mathcal{T}_{1}, \mathcal{T}_{3}$ be two disjoint sets of primes satisfying
$$ \mathcal{T}_{i} \subseteq [(\log^{(1/\lambda_{1})+\epsilon}X)/2,\log^{(1/\lambda_{1})+\epsilon}X] \;\;\; \textrm{ and } \;\;\; \#\mathcal{T}_{i} \geq \frac{\log^{(1/\lambda_{1})+\epsilon}X}{5 \log(\log^{(1/\lambda_{1})+\epsilon}X)}. $$
Then we can let $\mathcal{C} \subseteq [1,X]$, $\mathcal{B} \subseteq [1,Y]$ and $\mathcal{A} \subseteq [1,Z]$ be the subsets of squarefree numbers having all their prime factors from $\mathcal{T}_{1},\mathcal{T}_{2},\mathcal{T}_{3}$ respectively. By Lemma 2 these sets have size at least $X^{1-\lambda_{1}+\Omega(\epsilon)+o(1)}, Y^{1-\lambda_{1}+\Omega(\epsilon)+o(1)}, Z^{1-\lambda_{1}+\Omega(\epsilon)+o(1)}$, and one can check in the proof of Lemma 2 that this will still hold if we remove numbers smaller than $X^{1-\delta},Y^{1-\delta},Z^{1-\delta}$ (respectively) from the sets. Clearly each element of $\mathcal{C}$ is coprime to each element of $\mathcal{A}$, and for each $a \in \mathcal{A}$ we have
$$ \#\{b \in \mathcal{B} : (b+1,a) \neq 1\} \leq \sum_{p \in \mathcal{T}_{3}} \#\{b \in \mathcal{B} : b \equiv -1 \textrm{ mod } p\} \ll \sum_{p \in \mathcal{T}_{3}} \frac{\#\mathcal{B}}{p} \leq \#\mathcal{B}/2 $$
provided $X$ is sufficiently large, using any of several results about the distribution of smooth numbers in arithmetic progressions. See e.g. Granville's paper~\cite{granville}.

This completes our construction of sets $\mathcal{A}, \mathcal{B}, \mathcal{C}$ producing many solutions of $cw-au=b+1$, and all of whose elements have their prime factors from a set of size $\leq \log^{(1/\lambda_{1})+\epsilon}X$.

\subsection{An unconditional proof}
Now we shall prove Theorem 2 without the aid of any conjecture. We shall need the following result, which is essentially Lemma 3 of Friedlander and Iwaniec's paper~\cite{friediwan}:
\begin{nthres3}[Friedlander and Iwaniec, 1985]
For any $q \geq 2$, any $N \geq 1$, and any $\eta > 0$,
$$ \frac{1}{\phi(q)} \sum_{\chi \textrm{ mod } q, \atop \chi \neq \chi_{0}} \left|\sum_{1 \leq n \leq N} \chi(n) \right|^{4} \ll_{\eta} q^{\eta} N^{2}, $$
where $\chi_{0}$ denotes the principal Dirichlet character to modulus $q$, and the implicit constant depends on $\eta$ only.
\end{nthres3}

Looking to bound the ``big Oh'' term from the end of $\S 4.1$, we see that for each $a \in \mathcal{A}$,
\begin{eqnarray}
&& \frac{1}{\phi(a)} \sum_{\chi \textrm{ mod } a, \atop \chi \neq \chi_{0}} \left|\sum_{b \in \mathcal{B}} \overline{\chi}(b+1) \right| \left|\sum_{c \in \mathcal{C}} \chi(c) \right| \left|\sum_{1 \leq w \leq W} \chi(w) \right| \nonumber \\
& \leq & \sqrt{\frac{1}{\phi(a)} \sum_{\chi \textrm{ mod } a} \left|\sum_{b \in \mathcal{B}} \overline{\chi}(b+1) \right|^{2}} \left(\frac{1}{\phi(a)} \sum_{\chi \textrm{ mod } a} \left|\sum_{c \in \mathcal{C}} \chi(c) \right|^{4} \right)^{1/4} \left(\frac{1}{\phi(a)} \sum_{\chi \textrm{ mod } a, \atop \chi \neq \chi_{0}} \left|\sum_{1 \leq w \leq W} \chi(w) \right|^{4} \right)^{1/4} \nonumber \\
& \ll_{\eta} & Z^{\eta} \sqrt{W} \sqrt{\#\mathcal{B}+Y^{1+\eta}/Z^{1-\delta}} \left(\frac{1}{\phi(a)} \sum_{\chi \textrm{ mod } a} \left|\sum_{c \in \mathcal{C}} \chi(c) \right|^{4} \right)^{1/4}, \nonumber
\end{eqnarray}
using two applications of the Cauchy--Schwarz inequality, Number Theory Result 3, and the same argument as in $\S 4.2$ to estimate the sums involving $\mathcal{B}$. The sums involving $\mathcal{C}$ have size at most the number of quadruples $(c_{1},c_{2},c_{3},c_{4}) \in \mathcal{C}^{4}$ for which $a|(c_{1}c_{2}-c_{3}c_{4})$, and this is $\ll_{\eta} X^{\eta}((\#\mathcal{C})^{2}+X^{2}/Z^{1-\delta})$, exactly similarly to the argument in $\S 4.2$ (where the $X^{\eta}(\#\mathcal{C})^{2}$ term counts quadruples with $c_{1}c_{2}=c_{3}c_{4}$). Actually, since we will end up choosing $\mathcal{C}$ such that $\#\mathcal{C} \leq \sqrt{X}$, and we have been assuming since $\S 4.1$ that $Z \leq X$, the term $(\#\mathcal{C})^{2}$ is redundant in this bound.

\vspace{12pt}
Now, exactly similarly to $\S 4.2$, we need to choose $W,X,Y,Z$ and the sets $\mathcal{A},\mathcal{B},\mathcal{C}$ so that both
\begin{eqnarray}
&& \sum_{a \in \mathcal{A}} \frac{1}{\phi(a)} \#\{b \in \mathcal{B} : (b+1,a)=1\} \#\{c \in \mathcal{C} : (c,a) = 1\} \#\{1 \leq w \leq W : (w,a)=1\} \nonumber \\
& \geq & B(\delta) \sqrt{W} \#\mathcal{A} \sqrt{(\#\mathcal{B} + Y/Z)(X/\sqrt{Z})}, \nonumber
\end{eqnarray}
and also
\begin{eqnarray}
&& \sum_{a \in \mathcal{A}} \frac{1}{\phi(a)} \#\{b \in \mathcal{B} : (b+1,a)=1\} \#\{c \in \mathcal{C} : (c,a) = 1\} \#\{1 \leq w \leq W : (w,a)=1\} \nonumber \\
& \geq & \frac{B(\delta) XW^{2}}{Z}. \nonumber
\end{eqnarray}
In particular, if $\alpha \geq 0$ is such that  $\#\mathcal{A} \geq Z^{1-\alpha}$, $\mathcal{B} \geq Y^{1-\alpha}$ and $\#\mathcal{C} \geq X^{1-\alpha}$, as in $\S 4.2$, then {\em it will suffice if}
$$ \sqrt{W} Y^{1-\alpha} \geq B(\delta) X^{\alpha-1/2} Z^{3/4} \sqrt{\max\{Y^{1-\alpha},Y/Z\}} \;\;\; \textrm{ and } \;\;\; Y^{1-\alpha}Z^{1-\alpha} \geq B(\delta) X^{\alpha}W. $$

If we choose $W$ as large as possible, subject to the second constraint, then the first constraint becomes (after squaring both sides)
$$ (Y^{1-\alpha}Z^{1-\alpha}X^{-\alpha}) Y^{2(1-\alpha)} \geq B(\delta) X^{2\alpha-1} Z^{3/2} \max\{Y^{1-\alpha},Y/Z\}. $$
Comparing the powers of $Y$ on either side, since we shall certainly end up with $\alpha \leq 2/3$ the power on the left will be greater than on the right, and so we should choose $Y$ as large as we can. We have been supposing since $\S 4.1$ that $Y \leq B(\delta) XW$, so if we set $Y=B(\delta) XW=B(\delta) X^{1-\alpha}Y^{1-\alpha}Z^{1-\alpha}$, i.e. $Y^{\alpha} = B(\delta) X^{1-\alpha} Z^{1-\alpha}$, then our constraint becomes (after raising both sides to the power $\alpha$)
$$ Z^{\alpha(1-\alpha)} X^{-\alpha^{2}} (X^{1-\alpha}Z^{1-\alpha})^{3(1-\alpha)} \geq B(\delta) X^{2\alpha^{2}-\alpha} Z^{3\alpha/2} \max\{(X^{1-\alpha}Z^{1-\alpha})^{1-\alpha},(X^{1-\alpha}Z^{1-\alpha})/Z^{\alpha}\}. $$
If we shift all powers of $X$ and $Z$ to the right, and simplify as much as possible, this remaining constraint becomes
$$ \max\{X^{-2+3\alpha+\alpha^{2}}Z^{-2+9\alpha/2-\alpha^{2}}, X^{-2+4\alpha}Z^{-2+9\alpha/2-2\alpha^{2}}\} \leq B(\delta). $$
When $1/2 \leq \alpha \leq 1$ we have $-2+9\alpha/2-\alpha^{2} \geq 0$, so if the above is to hold we must have $X^{(-2+4\alpha)(-2+9\alpha/2-\alpha^{2})} \leq B(\delta) Z^{(2-9\alpha/2+2\alpha^{2})(-2+9\alpha/2-\alpha^{2})} \leq B(\delta) X^{(2-9\alpha/2+2\alpha^{2})(2-3\alpha-\alpha^{2})}$. Since $X$ is a large parameter, this amounts to asking that
$$ (-2+4\alpha)(-2+9\alpha/2-\alpha^{2}) < (2-9\alpha/2+2\alpha^{2})(2-3\alpha-\alpha^{2}), $$
or equivalently, after simplifying, that $-2\alpha + 9\alpha^{2}/2 - 5\alpha^{3}/2 + 2\alpha^{4} < 0$. Since we are looking to have $1/2 \leq \alpha \leq 1$ we can cancel one multiple of $\alpha$, and are left wanting $4\alpha^{3}-5\alpha^{2}+9\alpha-4 < 0$. This will be satisfied when $\alpha = \lambda_{0}-\epsilon$, by definition of $\lambda_{0}$ as the real root of that cubic.

If we choose $\alpha$ in this way; and put $Z = B(\delta) X^{(-2+4\alpha)/(2-9\alpha/2+2\alpha^{2})}$, and similarly choose $Y = B(\delta) X^{(-2\alpha^{2}+5\alpha/2-1/2)/(2-9\alpha/2+2\alpha^{2})}$ and $W = B(\delta) X^{(-4\alpha^{2}+7\alpha-5/2)/(2-9\alpha/2+2\alpha^{2})}$ as listed above; then our various inequalities and assumptions will all be satisfied. Theorem 2 follows on constructing the sets $\mathcal{A},\mathcal{B},\mathcal{C}$ exactly similarly to $\S 4.2$, simply replacing $1/\lambda_{1} + \epsilon$ by $1/\lambda_{0} + \epsilon$ in all of the exponents.
\begin{flushright}
Q.E.D.
\end{flushright}

\section{Further results and discussion}

\subsection{``Optimality'' of the arguments in $\S 3$}
The arguments in $\S 3$ established the existence of many solutions to the $S$-unit equation $a+1=c$ by counting solutions to an associated linear equation on average over its coefficients. To do this they made use of three kinds of tools, namely:
\begin{enumerate}
\item obvious uses of the orthogonality of characters $\chi$;

\item bounds for incomplete character sums $\sum_{n=M+1}^{M+N} \chi(n)$;

\item the multiplicative large sieve.
\end{enumerate}
The ways in which these were combined in the proofs may have seemed a bit arbitrary, but in this subsection we will sketch an argument that, roughly speaking, one cannot improve upon the results from $\S 3$ (or at least the conditional result in $\S 3.2$) if one only uses these tools. This suggests that if one wishes to make further progress towards the conjectural bound $e^{s^{1/2-\epsilon}}$ in Theorem 1, one will need either different general tools (some possibilities being discussed in the next subsection), or more specific information about e.g. character sums over smooth numbers. The reader will note that the only information we used about smooth numbers was an estimate for how many there are, as in Lemma 2, and the fact that we could factor certain sets of smooth numbers in useful ways.

\vspace{12pt}
We suppose we are arguing in the general framework of $\S 3$, so (as in $\S 3.2$ and $\S 3.3$) we are seeking the largest value of $\alpha$ for which both
$$ X^{1-\alpha}WZ^{-\alpha} \geq B(\delta) \sum_{a \in \mathcal{A}} \frac{1}{\phi(a)} \sum_{\chi \textrm{ mod } a, \atop \chi \neq \chi_{0}} \left|\sum_{c \in \mathcal{C}} \chi(c) \right| \left|\sum_{1 \leq w \leq W} \chi(w) \right| $$
and
$$ Z^{1-\alpha} X^{-\alpha} \geq B(\delta) W, $$
where $\mathcal{C} \subseteq [X^{1-\delta},X]$ and $\mathcal{A} \subseteq [Z^{1-\delta},Z]$ are sets satisfying $\#\mathcal{C} \geq X^{1-\alpha}$ and $\#\mathcal{A} \geq Z^{1-\alpha}$. We need to bound the character sums in the first inequality, and let us restrict ourselves to doing this by performing some combination of the following ``moves'':
\begin{itemize}
\item use the bound $\frac{1}{\phi(a)} \sum_{\chi \textrm{ mod } a} \left|\sum_{l \in \mathcal{L}} c_{l} \chi(l) \right|^{2} \leq \sum_{l \in \mathcal{L}} |c_{l}|^{2}$ for any set $\mathcal{L}$ contained in an interval of length $< a$, and any complex numbers $c_{l}$;

\item apply the P\'{o}lya--Vinogradov inequality and the large sieve, as in $\S 3.3$;

\item use the bound $\left|\sum_{1 \leq w \leq W} \chi(w) \right| \leq B(\delta) \sqrt{W}$, where $\chi$ is any non-principal character mod $a$.
\end{itemize}
The bound in the third move is stronger than we can presently prove, so assuming that we can freely use it will make our optimality result stronger. We will also assume that we can factor $\mathcal{C}$ as a product of sets $\mathcal{Q}_{1},...,\mathcal{Q}_{k}$, for any fixed $k$ that we wish, and a priori for any sizes of sets $\mathcal{Q}_{i}$ that we wish. (Here we mean that the elements of $\mathcal{C}$ should be the products $q_{1}q_{2}...q_{k}$, for $q_{i} \in \mathcal{Q}_{i}$, and that these products should be distinct).

Let $k \geq 2$, and let $\mathcal{I} \subseteq \{2,3,...,k\}$ be any subset. We think of the set $\mathcal{C}$ as having been factored into sets $\mathcal{Q}_{1},...,\mathcal{Q}_{k}$, and of $\mathcal{I}$ as denoting the set of indices $i$ for which we shall estimate the character sum over $\mathcal{Q}_{i}$ using a ``large sieve move''. These are the character sums that we shall need to pair up with a sum $|\sum_{1 \leq w \leq W} \chi(w)|$, so that we can use the P\'{o}lya--Vinogradov inequality and obtain the Gauss sum weights needed in the large sieve. (Recall $\S 3.3$). The character sums over the other $\mathcal{Q}_{j}$ will be estimated using an ``orthogonality move'', with the associated copy of $|\sum_{1 \leq w \leq W} \chi(w)|$ (which will arise from applying the Cauchy--Schwarz inequality) being estimated using our conjectural bound $B(\delta) \sqrt{W}$. Note that we may certainly assume, after possibly relabelling $\mathcal{Q}_{1}$ as $\mathcal{Q}_{2}$, that the large sieve is {\em not} used to estimate the character sum over $\mathcal{Q}_{1}$, because if it was we would end up with precisely the argument given in $\S 3.3$.

Now $\sum_{a \in \mathcal{A}} 1/\phi(a) \sum_{\chi \textrm{ mod } a, \atop \chi \neq \chi_{0}} |\sum_{c \in \mathcal{C}} \chi(c)| |\sum_{1 \leq w \leq W} \chi(w)|$ is
\begin{eqnarray}
&& \sum_{a \in \mathcal{A}} \frac{1}{\phi(a)} \sum_{\chi \textrm{ mod } a, \atop \chi \neq \chi_{0}} \left( \prod_{i=1}^{k} \left|\sum_{q \in \mathcal{Q}_{i}} \chi(q) \right| \right) \left|\sum_{1 \leq w \leq W} \chi(w) \right| \nonumber \\
& \leq & \sqrt{\sum_{a \in \mathcal{A}} \frac{1}{\phi(a)} \sum_{\chi \textrm{ mod } a, \atop \chi \neq \chi_{0}} \left|\sum_{q \in \mathcal{Q}_{1}} \chi(q) \right|^{2}} \cdot \prod_{i \in \mathcal{I}} \left( \sum_{a \in \mathcal{A}} \frac{1}{\phi(a)} \sum_{\chi \textrm{ mod } a, \atop \chi \neq \chi_{0}} \left|\sum_{q \in \mathcal{Q}_{i}} \chi(q) \right|^{2^{i}} \left|\sum_{1 \leq w \leq W} \chi(w) \right|^{2} \right)^{1/2^{i}} \cdot \nonumber \\
&& \cdot \prod_{2 \leq i \leq k, \atop i \notin \mathcal{I}} \left( \sum_{a \in \mathcal{A}} \frac{1}{\phi(a)} \sum_{\chi \textrm{ mod } a, \atop \chi \neq \chi_{0}} \left|\sum_{q \in \mathcal{Q}_{i}} \chi(q) \right|^{2^{i}} \right)^{1/2^{i}} \cdot \left( \sum_{a \in \mathcal{A}} \frac{1}{\phi(a)} \sum_{\chi \textrm{ mod } a, \atop \chi \neq \chi_{0}} \left|\sum_{1 \leq w \leq W} \chi(w) \right|^{2^{k} - \sum_{i \in \mathcal{I}} 2^{k-i+1} } \right)^{1/2^{k}}, \nonumber
\end{eqnarray}
by repeated applications of the Cauchy--Schwarz inequality. Given our limited set of moves, the only way we can handle the large powers of character sums is to note that
$$ \left|\sum_{q \in \mathcal{Q}_{i}} \chi(q) \right|^{2^{i}} = \left|\left( \sum_{q \in \mathcal{Q}_{i}} \chi(q) \right)^{2^{i-1}} \right|^{2} = \left| \sum_{Q_{i}^{2^{i-1}(1-\delta)} \leq l \leq Q_{i}^{2^{i-1}}} \#\{(l_{1},...,l_{2^{i-1}}) \in \mathcal{Q}_{i}^{2^{i-1}} : l_{1}...l_{2^{i-1}} = l\} \chi(l) \right|^{2} $$
if $\mathcal{Q}_{i} \subseteq [Q_{i}^{1-\delta},Q_{i}]$, say, where $\prod_{i=1}^{k} Q_{i} = X$. Thus if $i \notin \mathcal{I}$, so we will be using orthogonality, we need to have $Q_{i}^{2^{i-1}} < Z$, and then $1/\phi(a) \sum_{\chi \textrm{ mod } a} |\sum_{q \in \mathcal{Q}_{i}} \chi(q)|^{2^{i}}$ will be at most
\begin{eqnarray}
\sum_{l \leq Q_{i}^{2^{i-1}}} \#\{(l_{1},...,l_{2^{i-1}}) \in \mathcal{Q}_{i}^{2^{i-1}} : l_{1}...l_{2^{i-1}} = l\}^{2} & \ll_{\eta,i} & X^{\eta} \sum_{l} \#\{(l_{1},...,l_{2^{i-1}}) \in \mathcal{Q}_{i}^{2^{i-1}} : l_{1}...l_{2^{i-1}} = l\} \nonumber \\
& = & X^{\eta}(\#\mathcal{Q}_{i})^{2^{i-1}}, \nonumber
\end{eqnarray}
using a crude bound for the divisor function. Similarly if $i \in \mathcal{I}$, so we will be using the large sieve, we need to have $Q_{i}^{2^{i-1}} < Z^{2}$, and then we can have a bound
$$ \sum_{a \in \mathcal{A}} \frac{1}{\phi(a)} \sum_{\chi \textrm{ mod } a, \atop \chi \neq \chi_{0}} \left|\sum_{q \in \mathcal{Q}_{i}} \chi(q) \right|^{2^{i}} \left|\sum_{1 \leq w \leq W} \chi(w) \right|^{2} \ll_{\eta,i} X^{\eta} Z^{2} (\#\mathcal{Q}_{i})^{2^{i-1}}. $$

Applying these bounds and our conjectural estimate for $\sum_{1 \leq w \leq W} \chi(w)$, we see that our target inequalities become
$$ X^{1-\alpha}WZ^{-\alpha} \geq B(\delta) \sqrt{\prod_{i=1}^{k} (\#\mathcal{Q}_{i})} (\#\mathcal{A})^{1/2 + \sum_{2 \leq i \leq k, i \notin \mathcal{I}} 2^{-i} + 2^{-k}} Z^{2 \sum_{i \in \mathcal{I}} 2^{-i}} W^{1/2-\sum_{i \in \mathcal{I}} 2^{-i}}, $$
and 
$$ Z^{1-\alpha} X^{-\alpha} \geq B(\delta) W. $$
Here we have $\prod_{i=1}^{k} (\#\mathcal{Q}_{i}) = \#\mathcal{C} \geq X^{1-\alpha}$, since $\mathcal{Q}_{i}$ were the factors of $\mathcal{C}$, and we also have $X = \prod_{i=1}^{k} Q_{i} \approx Z^{1+\sum_{2 \leq i \leq k, i \notin \mathcal{I}} 2^{1-i} + \sum_{i \in \mathcal{I}} 2^{2-i}} = Z^{2-2^{1-k} + \sum_{i \in \mathcal{I}} 2^{1-i}} $ in view of our discussion about the permissible sizes of the $Q_{i}$. If we set $\theta = \sum_{i \in \mathcal{I}} 2^{-i}$, and simplify where possible, our target inequalities become
$$ X^{(1-\alpha)/2} W^{1/2 + \theta} \geq B(\delta) Z^{1+(1+\alpha)\theta} \;\;\; \textrm{ and } \;\;\; Z^{1-\alpha} X^{-\alpha} \geq B(\delta) W, $$
where $X \approx Z^{2-2^{1-k}+2\theta}$.

Now we can find the largest permissible value of $\alpha$ exactly as in all our previous calculations: it turns out that we must have
$$ \alpha < \frac{1-2^{1-k}+2\theta}{5-2^{2-k}+12\theta-2^{2-k}\theta+4\theta^{2}} = \frac{1}{5} \frac{(1-2^{1-k}+2\theta)}{(1-(2/5)2^{1-k}+(12/5)\theta-(1/5)2^{2-k}\theta+(4/5)\theta^{2})}. $$
In particular, the right hand side is $< 1/5$ whenever $k \geq 2$, so that we cannot use these methods to improve on the conditional bound $1/5-\epsilon$ from $\S 3.2$.

The reader might also note e.g. that for any given $k$, we must either have $\theta = 0$ or $\theta \geq 2^{-k}$, and so we always have $12\theta-2^{2-k}\theta+4\theta^{2} \geq 12\theta = 6 \cdot 2\theta$. This means that, {\em if} we have access to a strong bound $\left|\sum_{1 \leq w \leq W} \chi(w) \right| \leq B(\delta) \sqrt{W}$ for incomplete character sums, we will have the best chance of improving on our exponent $1/6-\epsilon$ (from $\S 3.3$) if we do not use the large sieve at all, so that $\theta = 0$.

\vspace{12pt}
The alert reader will note that we did not really need the bound $\left|\sum_{1 \leq w \leq W} \chi(w) \right| \leq B(\delta) \sqrt{W}$ in the above calculations, but rather a bound $B(\delta) W^{t}$ for the mean value $(1/\phi(a)) \sum_{\chi \textrm{ mod } a, \chi \neq \chi_{0}} \left|\sum_{1 \leq w \leq W} \chi(w) \right|^{2t}$, for certain natural numbers $t$. Arranging things so that we only needed a mean value estimate cost us one extra application of the Cauchy--Schwarz inequality, which explains why we did not recover the value $\alpha = 1/5 - \epsilon$ when $\theta = 0$. When $t=2$ such a mean value estimate is a known result, (indeed we stated it as Number Theory Result 3, in $\S 4.3$), but so far as the author is aware there is no known proof of such an estimate for any $t > 2$. Consequently we cannot rigorously implement the foregoing calculations for large values of $k$, except if the large sieve is used lots of times; and in that case $\theta$ will be quite large, which we just noted was undesirable for obtaining a strong result. An interested reader may care to explore for what value of $t$ one would need a sharp mean value bound to improve on our exponent $1/6 - \epsilon$.

\subsection{More results on counting solutions to linear equations on average}
In this subsection we will give two general results on solving linear equations on average over their coefficients. The author proved these whilst pursuing the work in $\S 3$, but as we shall discuss the results they supply about $S$-unit equations are not as strong as those already presented. We will also briefly discuss some other ideas for obtaining results of this kind, and mention some connections with work of other authors.

\vspace{12pt}
Unlike the arguments in $\S 3$, the following theorem can supply quite good information if one wants the coefficient sets $\mathcal{A}, \mathcal{C}$ to be on the same scale, i.e. if one wants to take $X$ and $Z$ about the same size:
\begin{thm3}
Suppose that $\mathcal{C} \subseteq [3X/4,X]$ and $\mathcal{A} \subseteq [3Z/4,Z]$ are any sets of {\em natural} numbers, and that the set $\mathcal{C}$ ``factors'' as a product of some sets $\mathcal{Q} \subseteq [3Q/4,Q]$ and $\mathcal{R} \subseteq [3R/4,R]$, in the sense that $QR=X$, and
$$ \mathcal{C} = \{qr : q \in \mathcal{Q}, r \in \mathcal{R}\}, $$
and all of these products $qr$ are distinct. Also suppose that $Z \leq X$, and that each element of $\mathcal{C}$ is coprime to each element of $\mathcal{A}$. Let $1/\sqrt{Z} \leq \mu \leq 1$ be any number, and suppose that $\mathcal{R}$ satisfies $\#\mathcal{R} \geq 1/\mu$. Then for any $\eta > 0$,
\begin{eqnarray}
\sum_{a \in \mathcal{A}} \sum_{1 \leq w \leq \mu a} \sum_{c \in \mathcal{C}} \textbf{1}_{cw-au=1 \textrm{ for some } u} & = & (1+O(1/\log Z))\mu \#\mathcal{A} \#\mathcal{C} + O(\#\mathcal{A} \sqrt{\#\mathcal{C} X (\log Z) / (\mu Z)}) + \nonumber \\
&& + O_{\eta}(X^{\eta} \sqrt{\mu \#\mathcal{A} \#\mathcal{C}} (\sqrt{QZ} + (XR)^{1/4}(ZR^{2}+Z^{2}+X/\mu)^{1/4})), \nonumber
\end{eqnarray}
where the constants implicit in the ``big Oh'' notation are absolute except in the third case, where the constant depends on $\eta$.
\end{thm3}

In earlier sections we summed over $1 \leq w \leq W$, whereas here the range of summation varies a little with $a$. This makes no real difference when applying Theorem 3, since all the elements $a \in \mathcal{A}$ are roughly the same size, but this formulation makes some steps of the proof a bit neater. The requirement that each element of $\mathcal{C}$ is coprime to each element of $\mathcal{A}$ could be removed, if the conclusion of the theorem was adjusted suitably.

The proof of Theorem 3 will occupy most of the rest of this paper. We apply the circle method with additive characters, rather than multiplicative characters, noting that the triple sum is
\begin{eqnarray}
\sum_{a \in \mathcal{A}} \sum_{1 \leq w \leq [\mu a]} \sum_{c \in \mathcal{C}} \frac{1}{a} \sum_{-a/2 < h \leq a/2} e(h(c^{-1}-w)/a) & = & \sum_{a \in \mathcal{A}} \frac{[\mu a] \#\mathcal{C}}{a} + \nonumber \\
&& + \sum_{a \in \mathcal{A}} \sum_{-a/2 < h \leq a/2, \atop h \neq 0} (\frac{1}{a} \sum_{1 \leq w \leq [\mu a]} e(\frac{-hw}{a})) \sum_{c \in \mathcal{C}} e(\frac{hc^{-1}}{a}). \nonumber
\end{eqnarray}
Here, as is usual, we let $e(t) = e^{2\pi i t}$, and used $c^{-1}$ to denote the multiplicative inverse of $c$ with respect to the relevant modulus (in this case $a$). Now for non-zero $-a/2 < h \leq a/2$,
$$ \frac{1}{a} \sum_{1 \leq w \leq [\mu a]} e(-hw/a) = \frac{1}{a} \frac{e(-h[\mu a]/a) - 1}{1 - e(h/a)} = \frac{1}{a} \frac{e(-h[\mu a]/a) - 1}{-2\pi i h/a} + O(1/a). $$
This is easily seen if one distinguishes the cases where $|h| \leq a/100$, say, and where $|h|$ is larger (and so the whole sum is $O(1)$). Moreover we have
$$ \frac{1}{a} \frac{e(-h[\mu a]/a) - 1}{-2\pi i h/a} + O(1/a) = \frac{1}{a} \frac{e(-h\mu) - 1}{-2\pi i h/a} + O(1/a) =: s_{\mu}(h) + O(1/a), $$
say, and so for any $\lambda < 3/8$ (so that $\lambda \mu Z < (1/2)\min_{a \in \mathcal{A}} a$) our triple sum is
\begin{eqnarray}
\mu \#\mathcal{A} \#\mathcal{C} & + & \sum_{a \in \mathcal{A}} \sum_{-\lambda \mu Z \leq h \leq \lambda \mu Z, \atop h \neq 0} s_{\mu}(h) \sum_{c \in \mathcal{C}} e(hc^{-1}/a) + O((\lambda \mu + 1/Z) \#\mathcal{A} \#\mathcal{C}) + \nonumber \\
& + & \sum_{a \in \mathcal{A}} \sum_{-a/2 < h \leq a/2, \atop |h| > \lambda \mu Z} (\frac{1}{a} \sum_{1 \leq w \leq [\mu a]} e(-hw/a)) \sum_{c \in \mathcal{C}} e(hc^{-1}/a). \nonumber
\end{eqnarray}

To deal with the sums over large $|h|$, we can fix $a \in \mathcal{A}$ and $H \leq Z$ and note that 
\begin{eqnarray}
\sum_{c \in \mathcal{C}} \left| \sum_{H < |h| \leq 2H} (\frac{1}{a} \sum_{1 \leq w \leq [\mu a]} e(\frac{-hw}{a})) e(\frac{hc^{-1}}{a}) \right| & \leq & \sqrt{\#\mathcal{C}} \sqrt{\sum_{c \in \mathcal{C}} \left| \sum_{H < |h| \leq 2H} (\frac{1}{a} \sum_{1 \leq w \leq [\mu a]} e(\frac{-hw}{a})) e(\frac{hc^{-1}}{a}) \right|^{2} } \nonumber \\
& \ll & \sqrt{\#\mathcal{C}} \sqrt{\frac{X}{a} \sum_{r=0}^{a-1} \left| \sum_{H < |h| \leq 2H} (\frac{1}{a} \sum_{1 \leq w \leq [\mu a]} e(\frac{-hw}{a})) e(\frac{hr}{a}) \right|^{2} }, \nonumber
\end{eqnarray}
on slicing the set $\mathcal{C}$ into $\Theta(X/a)$ pieces of length $a$. Then we can expand the square and perform the summation over $r$, obtaining a bound
$$ \ll \sqrt{\#\mathcal{C}} \sqrt{X \sum_{r=0}^{a-1} \frac{1}{H^{2}} \#\{ H < |h| \leq 2H : h \equiv r \textrm{ mod } a\}^{2} } \ll \sqrt{\#\mathcal{C} X / H} $$
since we have $H \leq Z \ll a$ and $|(1/a) \sum_{1 \leq w \leq [\mu a]} e(-hw/a)| \ll 1/H$. This estimate clearly implies that the contribution from $|h| > \lambda \mu Z$ to our triple sum is $ \ll \#\mathcal{A} \sqrt{\#\mathcal{C} X / (\lambda \mu Z)} $, and on choosing $\lambda = 1/\log Z$ this is acceptable\footnote{If this estimate was unacceptable for some application, one could probably reduce the contribution from large $|h|$ to a satisfactory level by counting the numbers $w$ with a smooth weight (rather than the sharp cutoff to the interval $1 \leq w \leq [\mu a]$). We leave it to the interested reader to work out the details.} in Theorem 3. The error term $O((\lambda \mu + 1/Z) \#\mathcal{A} \#\mathcal{C})$ is also $O((\mu \#\mathcal{A} \#\mathcal{C})/\log Z)$ with this choice, since we assumed that $\mu \geq 1/\sqrt{Z}$.

\vspace{12pt}
We are left with a neat looking sum of weighted Kloosterman sums, which we shall estimate using a deep result of Deshouillers and Iwaniec~\cite{desiw}. The author was inspired to do this by the very nice use of Kloosterman sum estimates in a paper by Matom\"{a}ki~\cite{mat} (concerned with a rather different problem). It also works out a suggestion in Question 32 of Shparlinski's survey paper~\cite{shpar2}. We will state the relevant result now, although it will require a little preparation before we can apply it.

\begin{nthres4}[Deshouillers and Iwaniec, 1982]
Let $C,D,N,R$ be numbers at least as large as $1/2$, and $\textbf{b}_{n,r}$ be a complex sequence. Let $g(c)$ be an infinitely continuously differentiable function, having compact support in $[C,2C]$, and satisfying
$$ \frac{d^{v}g}{d c^{v}} \ll_{v} c^{-v} $$
for any $v \geq 0$, the implicit constant depending on $v$ at most. Let $h(d)$, having compact support in $[D,2D]$, have the same properties. Then for any $\epsilon > 0$,
$$ \left| \sum_{R < r \leq 2R} \sum_{0 < n \leq N} \textbf{b}_{n,r} \sum_{c,d, \atop (rd,c)=1} g(c)h(d) e(n(rd)^{-1}/c) \right| \ll_{\epsilon} (CDNR)^{\epsilon} K(C,D,N,R) \sqrt{\sum_{R < r \leq 2R} \sum_{0 < n \leq N} |\textbf{b}_{n,r}|^{2} }, $$
where $(rd)^{-1}$ stands for a solution of $(rd)^{-1}rd \equiv 1$ (mod $c$), and
$$ K^{2}(C,D,N,R) = C(R+N)(C+DR) + C^{2}D\sqrt{(R+N)R} + D^{2}NR. $$
\end{nthres4}

This is a slight restatement of Theorem 12 of Deshouillers and Iwaniec~\cite{desiw}, where we have simplified the allowed type of smooth weight function somewhat, and have not mentioned an additional variable that one can have in their result\footnote{If one assumes that Selberg's eigenvalue conjecture is true this leads to an improvement of Number Theory Result 4. There seems to be a misprint in the statement of Theorem 12 of Deshouillers and Iwaniec~\cite{desiw}, but the author believes that on Selberg's conjecture one can replace $\sqrt{(R+N)R}$ by $\sqrt{NR}$ in the expression for $K^{2}(C,D,N,R)$. See below for a further remark about this.}.

Initially we shall work to bound $\sum_{a \in \mathcal{A}} \sum_{|h| \leq 1/\mu, h \neq 0} s_{\mu}(h) \sum_{c \in \mathcal{C}} e(hc^{-1}/a)$. Preparing to apply Number Theory Result 4, we note firstly that (very similarly to the proof of Lemma 11 of Matom\"{a}ki~\cite{mat}) 
\begin{eqnarray}
\left| \sum_{a \in \mathcal{A}} \sum_{|h| \leq \frac{1}{\mu}, \atop h \neq 0} s_{\mu}(h) \sum_{c \in \mathcal{C}} e(\frac{hc^{-1}}{a}) \right| & = & \left| \sum_{a \in \mathcal{A}} \sum_{q \in \mathcal{Q}} \sum_{|h| \leq \frac{1}{\mu}, \atop h \neq 0} s_{\mu}(h) \sum_{r \in \mathcal{R}} e(\frac{hq^{-1}r^{-1}}{a}) \right| \nonumber \\
& \leq & \sqrt{\#\mathcal{A} \#\mathcal{Q}} \cdot \nonumber \\
&& \cdot \sqrt{ \sum_{\frac{2Z}{3} \leq a \leq \frac{4Z}{3}} g_{Z}(a) \sum_{\frac{2Q}{3} \leq q \leq \frac{4Q}{3}, \atop (q,a)=1} h_{Q}(q) \left| \sum_{|h| \leq \frac{1}{\mu}, \atop h \neq 0} s_{\mu}(h) \sum_{r \in \mathcal{R}, \atop (r,a)=1} e(\frac{hq^{-1}r^{-1}}{a}) \right|^{2} }. \nonumber
\end{eqnarray}
Here $\textbf{1}_{[3Z/4,Z]} \leq g_{Z} \leq \textbf{1}_{[2Z/3,4Z/3]}$ and $\textbf{1}_{[3Q/4,Q]} \leq h_{Q} \leq \textbf{1}_{[2Q/3,4Q/3]}$ can be any compactly supported functions having the smoothness properties required in Number Theory Result 4. Examining the contents of the second squareroot sign, if we expand the square and reorganise the terms we find this is
\begin{eqnarray}
&& \sum_{|h_{1}|,|h_{2}| \leq \frac{1}{\mu}, \atop h_{1},h_{2} \neq 0} \sum_{r_{1}, r_{2} \in \mathcal{R}} s_{\mu}(h_{1}) \overline{s_{\mu}(h_{2})} \sum_{\frac{2Q}{3} \leq q \leq \frac{4Q}{3}} h_{Q}(q) \sum_{\frac{2Z}{3} \leq a \leq \frac{4Z}{3}, \atop (a,qr_{1}r_{2})=1} g_{Z}(a) e(\frac{q^{-1}r_{1}^{-1}r_{2}^{-1}(h_{1}r_{2}-h_{2}r_{1})}{a}) \nonumber \\
& = & \sum_{\frac{R^{2}}{2} \leq r \leq R^{2}} \sum_{\frac{-2R}{\mu} \leq n \leq \frac{2R}{\mu}, \atop n \neq 0} \left(\sum_{|h_{1}|,|h_{2}| \leq 1/\mu, h_{1},h_{2} \neq 0 \atop {r_{1}, r_{2} \in \mathcal{R}, r_{1}r_{2}=r, \atop h_{1}r_{2}-h_{2}r_{1}=n}} s_{\mu}(h_{1}) \overline{s_{\mu}(h_{2})} \right) \sum_{\frac{2Q}{3} \leq q \leq \frac{4Q}{3}} h_{Q}(q) \sum_{\frac{2Z}{3} \leq a \leq \frac{4Z}{3}, \atop (a,qr)=1} g_{Z}(a) e(\frac{q^{-1}r^{-1}n}{a})  \nonumber \\
&& + \sum_{|h_{1}|,|h_{2}| \leq \frac{1}{\mu}, \atop h_{1},h_{2} \neq 0} \sum_{r_{1}, r_{2} \in \mathcal{R}} s_{\mu}(h_{1}) \overline{s_{\mu}(h_{2})} \textbf{1}_{h_{1}r_{2}=h_{2}r_{1}} \sum_{\frac{2Q}{3} \leq q \leq \frac{4Q}{3}} h_{Q}(q) \sum_{\frac{2Z}{3} \leq a \leq \frac{4Z}{3}, \atop (a,qr_{1}r_{2})=1} g_{Z}(a). \nonumber
\end{eqnarray}
If we let $b(n,r)$ denote the bracketed term on the second line, and note that $|s_{\mu}(h)| \ll \mu$ for all $h$, we find the above is
$$ \sum_{\frac{R^{2}}{2} \leq r \leq R^{2}} \sum_{\frac{-2R}{\mu} \leq n \leq \frac{2R}{\mu}, \atop n \neq 0} b(n,r) \sum_{\frac{2Q}{3} \leq q \leq \frac{4Q}{3}} h_{Q}(q) \sum_{\frac{2Z}{3} \leq a \leq \frac{4Z}{3}, \atop (a,qr)=1} g_{Z}(a) e(\frac{q^{-1}r^{-1}n}{a}) + O_{\eta}(\mu \#\mathcal{R} X^{\eta} QZ), $$
where the term $X^{\eta}$ (for arbitrarily small $\eta > 0$) arises as a crude bound for the divisor function evaluated at $h_{1}r_{2}$.

By Number Theory Result 4, with the quantities $C,D,N,R$ set as $2Z/3, 2Q/3, 2R/\mu, R^{2}/2$, the remaining sums are
$$ \ll_{\epsilon} (XZ)^{\epsilon} \sqrt{Z(R^{2}+R/\mu)(Z+QR^{2}) + Z^{2}Q\sqrt{(R^{2}+R/\mu)R^{2}} + Q^{2}R^{3}/\mu} \sqrt{\sum_{n,r} |b(n,r)|^{2}}. $$
Moreover we have, using the Cauchy--Schwarz inequality,
\begin{eqnarray}
\sum_{n,r} |b(n,r)|^{2} & \leq & \sum_{\frac{-2R}{\mu} \leq n \leq \frac{2R}{\mu}, \atop n \neq 0} \sum_{\frac{R^{2}}{2} \leq r \leq R^{2}} \left(\sum_{r_{1}, r_{2} \in \mathcal{R}, r_{1}r_{2}=r} 1 \right) \left(\sum_{r_{1}, r_{2} \in \mathcal{R}, r_{1}r_{2}=r} \left| \sum_{|h_{1}|,|h_{2}| \leq 1/\mu, \atop {h_{1},h_{2} \neq 0 \atop h_{1}r_{2}-h_{2}r_{1}=n}} s_{\mu}(h_{1}) \overline{s_{\mu}(h_{2})} \right|^{2} \right) \nonumber \\
& \ll_{\eta} & X^{\eta} \mu^{4} \sum_{\frac{-2R}{\mu} \leq n \leq \frac{2R}{\mu}, \atop n \neq 0} \sum_{r_{1},r_{2} \in \mathcal{R}} \left( \sum_{|h_{1}|,|h_{2}| \leq 1/\mu, \atop {h_{1},h_{2} \neq 0 \atop h_{1}r_{2}-h_{2}r_{1}=n}} 1 \right)^{2} \nonumber \\
& \leq & X^{\eta} \mu^{4} \sum_{r_{1},r_{2} \in \mathcal{R}} \sum_{-1/\mu \leq h_{1},h_{2},h_{3},h_{4} \leq 1/\mu, \atop h_{1},h_{2},h_{3},h_{4} \neq 0} \textbf{1}_{r_{2}(h_{1}-h_{3})=r_{1}(h_{2}-h_{4})} \nonumber \\
& \ll_{\eta} & X^{\eta} \mu^{4} ((\#\mathcal{R}/\mu)^{2} + X^{\eta}\#\mathcal{R}/\mu^{3}). \nonumber
\end{eqnarray}
The new $X^{\eta}$ term in the last line arises as a crude bound for the number of divisors of $r_{2}(h_{1}-h_{3})$, when $h_{1} \neq h_{3}$. Finally, since $QR=X$, and by assumption $R \geq \#\mathcal{R} \geq 1/\mu$, our bound from Number Theory Result 4 simplifies to
\begin{eqnarray}
& \ll_{\epsilon} & (XZ)^{\epsilon} \sqrt{ZR^{2}(Z+XR) + Z^{2}XR + X^{2}R/\mu} \sqrt{\mu^{4} (\#\mathcal{R}/\mu)^{2}} \nonumber \\
& \ll_{\epsilon} & (XZ)^{\epsilon} \sqrt{XR} \sqrt{ZR^{2} + Z^{2} + X/\mu} (\mu \#\mathcal{R}) \nonumber
\end{eqnarray}
(where we noted that $Z^{2}R^{2} \leq Z^{2}QR^{2} = Z^{2}XR$). Checking back over our calculations, recalling that we need to take a squareroot and multiply by $\sqrt{\#\mathcal{A} \#\mathcal{Q}}$ (and replacing $(XZ)^{\epsilon}$ by $X^{\eta}$), we have shown that
\begin{eqnarray}
\left| \sum_{a \in \mathcal{A}} \sum_{|h| \leq \frac{1}{\mu}, \atop h \neq 0} s_{\mu}(h) \sum_{c \in \mathcal{C}} e(\frac{hc^{-1}}{a}) \right| & \ll_{\eta} & X^{\eta} \sqrt{\#\mathcal{A} \#\mathcal{Q}} \sqrt{\mu \#\mathcal{R} QZ + \sqrt{XR} \sqrt{ZR^{2} + Z^{2} + X/\mu} (\mu \#\mathcal{R})} \nonumber \\
& \leq & X^{\eta} \sqrt{\mu \#\mathcal{A} \#\mathcal{C}} (\sqrt{QZ} + (XR)^{1/4}(ZR^{2}+Z^{2}+X/\mu)^{1/4}), \nonumber
\end{eqnarray}
which is acceptable for Theorem 3. In the second line we noted that $\#\mathcal{C} = \#\mathcal{Q} \#\mathcal{R}$.

It only remains to bound $\sum_{a \in \mathcal{A}} \sum_{1/\mu < |h| \leq \mu Z/\log Z} s_{\mu}(h) \sum_{c \in \mathcal{C}} e(hc^{-1}/a)$, presuming that the sum over $h$ is non-empty. However, the reader may check that if we split this into sums over dyadic ranges $H < |h| \leq 2H$, beginning with $H = 1/\mu$, and repeat the above argument on each range, then we always obtain at least as strong a bound as we did for the sum over $0 < |h| \leq 1/\mu$. Indeed, the only thing that changes is that we must sum over $|n| \leq 4RH$ rather than over $|n| \leq 2R/\mu$, and that we have the stronger bound $|s_{\mu}(h)| \ll 1/H$ rather than $|s_{\mu}(h)| \ll \mu$, and the latter always at least compensates for the former. Theorem 3 then follows.
\begin{flushright}
Q.E.D.
\end{flushright}

Unfortunately Theorem 3 is not a completely general result, since it requires one to be able to factor the coefficient set $\mathcal{C}$. This may seem a slightly bizarre state of affairs, but in fact one encounters related issues in many different settings: for example, Bombieri, Friedlander and Iwaniec~\cite{bfi1,bfi2,bfi3} were able to extend the range of summation over moduli in the famous Bombieri--Vinogradov theorem, but only in the presence of certain ``well-factorable'' weights that made it possible to apply results about averages of Kloosterman sums. There is related work, requiring the same well-factorable weights, on several important problems of analytic number theory.

\vspace{12pt}
Before moving on, we mention a few other ways in which one might proceed if studying linear equations using additive characters. Shparlinski~\cite{shpar} applies a deep estimate for bilinear forms with Kloosterman fractions due to Duke, Friedlander, and Iwaniec~\cite{dfi}, and thereby manages to obtain some information even if the sets $\mathcal{A}, \mathcal{C}$ are on the same scale and do not factor. However, the saving over trivial bounds using that result is something like $(X+Z)^{1/48}$, which is much less good than in Theorem 3. If one has sets that do factor, one could also try applying the Kloosterman sum estimates of Karatsuba~\cite{karatsuba}, which have the advantage that they do not require one to average over $a \in \mathcal{A}$. (Or, from another point of view, the disadvantage that they can't exploit any such averaging). In fact there are a great variety of estimates available for different kinds of averages of Kloosterman sums, many of which the author has not explored but may be of use in this context.

We remark that using Theorem 3, the author was able to prove Theorem 1 with a lower bound $e^{s^{1/10-\epsilon}}$, by setting $Z=Q=X^{2/3}$, $R=X^{1/3}$, and $\mu = Z^{-1/4} = X^{-1/6}$. Slightly remarkably, if one works on these scales then it makes no difference to the quality of this result whether one uses Theorem 3 or a version conditional on Selberg's eigenvalue conjecture (which would slightly improve the estimate of Number Theory Result 4). Using Theorem 3 of Karatsuba~\cite{karatsuba}, the author believes that our optimal parameter choices for Theorem 1 are $Z=R=Q=X^{1/2}$ and $\mu = Z^{-3/10} = X^{-3/20}$, and this results in a bound no better than $e^{s^{1/10-\epsilon}}$.

\vspace{12pt}
We finish by stating the following result, which can supply useful information when the sets $\mathcal{A}, \mathcal{C}$ are on different scales but need not factor:
\begin{thm4}
Suppose that $\mathcal{C} \subseteq [X/2,X]$ and $\mathcal{A} \subseteq [Z/2,Z]$ are any sets of squarefree natural numbers, and that $W \leq Z$ is any natural number. Then for any $\eta > 0$,
\begin{eqnarray}
\sum_{a \in \mathcal{A}} \sum_{1 \leq w \leq W} \sum_{c \in \mathcal{C}} \textbf{1}_{cw-au=1 \textrm{ for some } u} & = & \sum_{a \in \mathcal{A}} \frac{1}{\phi(a)} \#\{c \in \mathcal{C} : (c,a) = 1\} \#\{1 \leq w \leq W : (w,a)=1\} + \nonumber \\
&& + O_{\eta}((XZ)^{\eta} \sqrt{(\#\mathcal{A}) (\#\mathcal{C}) W \max\{X/Z, Z\}}), \nonumber
\end{eqnarray}
where the constant implicit in the ``big Oh'' notation depends on $\eta$ only.
\end{thm4}

The proof is similar to that from $\S 3.3$, but with one additional ingredient in bounding the error term. Thus we have 
\begin{eqnarray}
\sum_{a \in \mathcal{A}} \frac{1}{\phi(a)} \sum_{\chi \textrm{ mod } a, \atop \chi \neq \chi_{0}} \left|\sum_{c \in \mathcal{C}} \chi(c) \right| \left|\sum_{1 \leq w \leq W} \chi(w) \right| & \leq & \sqrt{\sum_{a \in \mathcal{A}} \frac{1}{\phi(a)} \sum_{\chi \textrm{ mod } a, \atop \chi \neq \chi_{0}} \frac{|\sum_{1 \leq w \leq W} \chi(w)|^{2}}{|\tau(\chi)|^{2}}} \cdot \nonumber \\
&& \cdot \sqrt{\sum_{a \in \mathcal{A}} \frac{1}{\phi(a)} \sum_{\chi \textrm{ mod } a, \atop \chi \neq \chi_{0}} \left|\sum_{c \in \mathcal{C}} \chi(c) \right|^{2} |\tau(\chi)|^{2}} \nonumber \\
& \ll_{\eta} & \sqrt{\sum_{a \in \mathcal{A}} \frac{1}{\phi(a)} \sum_{\chi \textrm{ mod } a, \atop \chi \neq \chi_{0}} \frac{|\sum_{1 \leq w \leq W} \chi(w)|^{2}}{|\tau(\chi)|^{2}}} \cdot \nonumber \\
&& \cdot (XZ)^{\eta} \sqrt{\max\{X,Z^{2}\} \#\mathcal{C}}, \nonumber
\end{eqnarray}
using the Cauchy--Schwarz inequality and Number Theory Result 2, and noting that $\tau(\chi)$ never vanishes if $\chi$ is a character to squarefree modulus. Seeking a non-trivial estimate for the sums under the first squareroot, we note they are
\begin{eqnarray}
\sum_{2 \leq y \leq Z} \sum_{\chi^{*} \textrm{ mod } y, \atop \chi^{*} \textrm{ primitive}} \sum_{a \in \mathcal{A}} \frac{1}{\phi(a)} \sum_{\chi \textrm{ mod } a, \atop \chi \textrm{ induced by } \chi^{*}} \frac{|\sum_{1 \leq w \leq W} \chi(w)|^{2}}{y} & \ll_{\eta} & Z^{\eta} \sum_{2 \leq y \leq Z} \sum_{\chi^{*} \textrm{ mod } y, \atop \chi^{*} \textrm{ primitive}} \sum_{a \in \mathcal{A}} \frac{1}{\phi(a)} \cdot \nonumber \\
&& \cdot \sum_{\chi \textrm{ mod } a, \atop \chi \textrm{ induced by } \chi^{*}} \frac{\sum_{d|a} |\sum_{1 \leq v \leq W/d} \chi^{*}(v)|^{2}}{y}, \nonumber
\end{eqnarray}
on employing the simple sieve-type identity
$$ \sum_{1 \leq w \leq W} \chi(w) = \sum_{1 \leq w \leq W, (w,a)=1} \chi^{*}(w) = \sum_{1 \leq w \leq W} \chi^{*}(w) \sum_{d|(w,a)} \mu(d) = \sum_{d | a} \mu(d) \sum_{1 \leq v \leq W/d} \chi^{*}(vd) $$
(where $\chi$ is a character mod $a$, induced by $\chi^{*}$, and $\mu$ denotes the M\"{o}bius function).

If $y$ divides $a$ then any primitive character mod $y$ will induce precisely one character mod $a$, whilst if $y$ does not divide $a$ it clearly cannot induce any. Thus the above is, except for a $Z^{\eta}$ multiplier,
\begin{eqnarray}
&& \sum_{2 \leq y \leq Z} \sum_{\chi^{*} \textrm{ mod } y, \atop \chi^{*} \textrm{ primitive}} \sum_{a \in \mathcal{A}} \frac{1}{\phi(a)} \frac{\sum_{d|a} |\sum_{1 \leq v \leq W/d} \chi^{*}(v)|^{2}}{y} \textbf{1}_{y|a} \nonumber \\
& = & \sum_{a \in \mathcal{A}} \frac{1}{\phi(a)} \sum_{d|a} \left(\sum_{2 \leq y \leq W, \atop y|a} \frac{1}{y} \sum_{\chi^{*} \textrm{ mod } y, \atop \chi^{*} \textrm{ primitive}} |\sum_{1 \leq v \leq W/d} \chi^{*}(v)|^{2} + \sum_{W < y \leq Z, \atop y|a} \frac{1}{y} \sum_{\chi^{*} \textrm{ mod } y, \atop \chi^{*} \textrm{ primitive}} |\sum_{1 \leq v \leq W/d} \chi^{*}(v)|^{2} \right) \nonumber \\
& \leq & \sum_{a \in \mathcal{A}} \frac{1}{\phi(a)} \sum_{d|a} \left(\sum_{2 \leq y \leq W, \atop y|a} \phi(y) \log^{2}y + \sum_{W < y \leq Z, \atop y|a} \frac{1}{y} \sum_{\chi^{*} \textrm{ mod } y, \atop \chi^{*} \textrm{ primitive}} |\sum_{1 \leq v \leq W/d} \chi^{*}(v)|^{2} \right), \nonumber
\end{eqnarray}
using Number Theory Result 1. Finally if we extend the summations over primitive characters to summations over all characters mod $y$, and use orthogonality, and apply a crude bound for the divisor function twice, we see the above is at most
\begin{eqnarray}
\sum_{a \in \mathcal{A}} \frac{1}{\phi(a)} \sum_{d|a} \left(\sum_{2 \leq y \leq W, \atop y|a} \phi(y) \log^{2}y + \sum_{W < y \leq Z, \atop y|a} \frac{1}{y} \phi(y) \frac{W}{d} \right) & \leq & \sum_{a \in \mathcal{A}} \frac{W}{\phi(a)} \sum_{d|a} \left(\sum_{2 \leq y \leq W, \atop y|a} \log^{2}y + \sum_{W < y \leq Z, \atop y|a} 1 \right) \nonumber \\
& \ll_{\eta} & Z^{\eta} \sum_{a \in \mathcal{A}} \frac{W}{\phi(a)}. \nonumber
\end{eqnarray}
Since this is $\ll_{\eta} Z^{\eta} (\#\mathcal{A}) W/Z$, Theorem 4 follows.

\begin{flushright}
Q.E.D.
\end{flushright}

The author did not try to incorporate the refinement in the above proof into the arguments of $\S 3$, and it is possible that it could lead to an improvement of those results. However, one could not directly combine it with the splitting of the set $\mathcal{C}$ in those arguments, (presumably one would have to use the Cauchy--Schwarz inequality again to separate $|\sum_{1 \leq w \leq W} \chi(w)|/|\tau(\chi)|$ from everything else), and the author suspects that one would lose more in doing this than one would gain.

The only obvious alternative approach to obtaining a result like Theorem 4 is to apply the additive circle method, in which case one would need to bound
$$ \sum_{a \in \mathcal{A}} \sum_{|h| \leq 1/\mu, \atop h \neq 0} s_{\mu}(h) \sum_{c \in \mathcal{C}} e(hc^{-1}/a) $$
along with various other terms. As suggested in Shparlinski's paper~\cite{shpar}, one could do this using Theorem 1 of Duke, Friedlander and Iwaniec~\cite{dfi}, and indeed it is not obvious what other bound could be applied if neither of the sets $\mathcal{A}, \mathcal{C}$ are assumed to factor. However, if one applies that theorem one obtains that
\begin{eqnarray}
\left| \sum_{a \in \mathcal{A}} \sum_{|h| \leq 1/\mu, \atop h \neq 0} s_{\mu}(h) \sum_{c \in \mathcal{C}} e(hc^{-1}/a) \right| & \ll & \mu \sum_{|h| \leq 1/\mu, \atop h \neq 0} \left| \sum_{a \in \mathcal{A}} \sum_{c \in \mathcal{C}} e(hc^{-1}/a) \right| \nonumber \\
& \ll_{\eta} & (XZ)^{\eta} \mu \frac{1}{\mu} \sqrt{\#\mathcal{A} \#\mathcal{C}}(\sqrt{X+Z} + \min\{X,Z\}), \nonumber
\end{eqnarray}
and this is never better than the error bound in Theorem 4 on the whole range $1 \leq Z \leq X$ (and is appreciably poorer in the interesting case where $W$ is small).

\vspace{12pt}
\noindent {\em Acknowledgements.} The author would like to thank his PhD supervisor, Ben Green, for his interest and encouragement.

\end{document}